\documentclass[11pt]{amsart}

\usepackage{amssymb,amsmath,latexsym}
\usepackage[dvips]{color}
\usepackage{xcolor}
\usepackage{comment}
\allowdisplaybreaks

\begin{document}
\newtheorem{thm}{Theorem}[section]
\newtheorem{lemma}[thm]{Lemma}
\newtheorem{defn}{Definition}[section]
\newtheorem{prop}[thm]{Proposition}
\newtheorem{corollary}[thm]{Corollary}
\newtheorem{remark}[thm]{Remark}
\newtheorem{example}[thm]{Example}
\newtheorem{assumption}{Assumption}
\numberwithin{equation}{section}

\def\ee{\varepsilon}
\def\qed{{\hfill $\Box$ \bigskip}}
\def\MM{{\cal M}}
\def\BB{{\cal B}}
\def\LL{{\cal L}}
\def\FF{{\cal F}}
\def\EE{{\cal E}}
\def\QQ{{\cal Q}}
\def\AA{{\cal A}}

\def\cB{\mbox{${\cal B}$}}
\def\<{\langle}  \def\>{\rangle}

\def\R{{\mathbb R}}
\def\N{{\mathbb N}}
\def\E{{\bf E}}
\def\F{{\bf F}}
\def\H{{\bf H}}
\def\P{{\bf P}}
\def\Q{{\bf Q}}
\def\S{{\bf S}}
\def\J{{\bf J}}
\def\K{{\bf K}}
\def\F{{\bf F}}
\def\A{{\bf A}}
\def\loc{{\bf loc}}
\def\eps{\varepsilon}
\def\semi{{\bf semi}}
\def\wh{\widehat}
\def\pf{\noindent{\bf Proof} }
\def\dim{{\rm dim}}

\title{Convergence rate for a class of supercritical superprocesses}
\author[R. Liu]{Rongli Liu}
\thanks{The research of R. Liu is supported in part by NSFC (Grant No. 11301261), and the Fundamental Research Funds for the Central Universities (Grant No.  2017RC007).}
\author[Y.-X. Ren]{Yan-Xia Ren}
\thanks{The research of Y.-X. Ren is supported in part  by NSFC (Grant Nos. 11731009 and 12071011),
 the National Key R\&D Program of China (No. 2020YFA0712900), and LMEQF}
\author[R. Song]{Renming Song}\thanks{The research of R. Song is supported in part by a grant from the Simons Foundation (\#429343, Renming Song).}
 \date{}\maketitle

\centerline{\bf Abstract}
Suppose $X=\{X_t, t\ge 0\}$ is a supercritical superprocess. Let $\phi$ be the non-negative eigenfunction of the mean semigroup of $X$ corresponding to the principal eigenvalue $\lambda>0$. Then $M_t(\phi)=e^{-\lambda t}\langle\phi, X_t\rangle, t\geq 0,$ is a non-negative martingale with almost sure limit $M_\infty(\phi)$.  In this paper we study the rate at which $M_t(\phi)-M_\infty(\phi)$ converges to $0$ as $t\to \infty$ when the process may not have finite variance. Under some conditions on the mean semigroup, we provide sufficient and  necessary conditions for the rate
in the almost sure sense. Some results on the convergence rate in $L^p$ with $p\in(1, 2)$ are also obtained.

\bigskip
\noindent{\bf Keywords}\hspace{2mm}
supercritical superprocess, convergence rate, infinite variance, spine decomposition,
 principal eigenvalue, eigenfunction, martingale
\medskip

\noindent{\bf  2010 MR Subject Classification}
primary 60J68; secondary 60F15; 60F25; 60G57

\bigskip
\begin{doublespace}
\section{Introduction and Main Results}

Let $\{Z_n, n\ge 0\}$ be a Galton-Watson process with $Z_0=1$ and  offspring mean $m:=EZ_1>1$, and let $W_n:=m^{-n}Z_n$. Then $\{W_n, n\ge 0\}$ is a non-negative martingale with almost sure limit $W_\infty$. It is well-known that $W_n$ converges to $W_\infty$ in $L^1$ if and only if
$E(Z_1\log^+ Z_1)<\infty$.
In the case $E(Z_1\log^+ Z_1)<\infty$,
it is natural to consider the rate at which $W_\infty-W_n$ converges to $0$.
In this paper we are mainly concerned with the convergence rate in
the almost sure sense,
when the process may not have finite variance.
This type of results first appeared in  Asmussen  \cite{A}, and then  in the book of Asmussen and Hering \cite{AH}.  The following result is from
\cite[Theorem II.4.1, p. 36]{AH}:
{\bf Theorem A.}
(i) Let  $p\in(1, 2)$ and $1/p + 1/q = 1$. Then
\begin{equation}\label{exp-rate}W_\infty -W_n = o(m^{-n/q})\quad\mbox{ a.s. as } n\to\infty\end{equation}
if and only if
$ E(Z^p_1)<\infty.$

(ii)
Let $\alpha>0$. Then $$\sum^\infty_{n=1}n^{\alpha-1}(W_\infty-W_n)\quad\mbox{ converges a.s.}$$
if and only if
$ E(Z_1(\log^+Z_1)^{1+\alpha})<\infty.$

(iii) Let $\alpha>0$. Then
$W_\infty-W_n=o(n^{-\alpha})$ a.s. as $n\to\infty$
if and only if
$$
E\Big[Z_1\Big(\log Z_1-\log n\Big)I_{\{Z_1>n\}}\Big]=
o([\log n]^{-\alpha}), \quad \mbox{as } n\to\infty.
$$

Asmussen  \cite{A} also discussed corresponding results for finite type  Galton-Watson processes, and continuous time Galton-Watson processes.
For multigroup branching diffusions on bounded domains,
convergence rate corresponding to Theorem A (i) is considered
in \cite[Section 13, Chapter VIII]{AH}. A sufficient condition, corresponding to $E(Z_1^p)<\infty$, is given for \eqref{exp-rate} to hold, see
\cite[Theorem VIII.13.2, p.343]{AH}.
The goal of this paper is to prove  the counterparts of the results in  Theorem A
for a class of superprocesses.

Before we give our model and results, we first review some
related work in the literature.
For any $p>1$, the $L^p$ convergence rate of $W_n-W_\infty$ to $0$ is obtained in Liu \cite[Proposition 1.3]{Liu01}. Huang and Liu \cite{HL} obtained  $L^p$ convergence rates for similar martingales in quenched and annealed senses for branching processes in random environment.
In \cite{AIPR,I,IM},  a class of non-negative intrinsic martingales $W_n$ for supercritical branching random walks were investigated. Let $W_\infty$ be the
almost sure limit of $W_n$ as $n\to\infty$.
Necessary and sufficient conditions for the $L^p$-convergence, $p>1$, of the series
\begin{equation}\label{exp-rate-BRW}
\sum^\infty_{n=1}e^{an}(W_\infty-W_n), \quad a>0,
\end{equation}
were obtained in \cite{AIPR},
 which may be viewed as the exponential rate of convergence of $E|W_\infty-W_n|^p$ to 0 as $n\to\infty$. In \cite{I}, sufficient conditions for the almost sure convergence of the series
\begin{eqnarray}\label{id:discrete}
\sum_{n=0}^\infty f(n)(W_\infty-W_n)
\end{eqnarray}
were obtained, where $f$ is a function regularly varying at $\infty$ with index larger than $-1$.
\cite{IM}  investigated sufficient conditions for \eqref{exp-rate-BRW} to
converge in the almost sure sense.
For general supercritical indecomposable multi-type branching processes,  sufficient conditions for polynomial rate of convergence in the
sense of convergence in probability
were given in \cite{IM15}.

We now introduce the setup of this paper.
We always assume that $E$ is a locally compact separable metric space.
 We will use $E_{\partial}:=E\cup\{\partial \}$ to denote
the one-point compactification of $E$. We will use ${\mathcal B}(E)$ and ${\mathcal B}(E_{\partial})$ to denote the Borel $\sigma$-fields on $E$
and $E_{\partial}$ respectively. ${\mathcal B}_b(E)$ (respectively ${\mathcal B}^+(E)$, respectively ${\mathcal B}^+_b(E)$) will denote the set of all bounded
(respectively non-negative, respectively bounded and non-negative) real-valued
Borel functions on $E$. All functions $f$ on $E$ will be automatically extended to $E_{\partial}$ by setting
$f(\partial)=0$.

We will always assume that
$\xi=\{(\xi_t)_{t\geq 0}; \Pi_x, x\in E\}$
is a Hunt process on $E$ and
$\zeta=\inf\{t>0:\,  \xi_t=\partial\}$ is the lifetime of $\xi$.
We use $(P_t)_{t\geq 0}$ to denote the semigroup of $\xi$ acting on functions defined on $E$ and $(\overline{P}_t)_{t\geq 0}$ to denote the semigroup of $\xi$ acting on functions defined on $E_\partial$.
We mention in passing that it is important that we take $E_{\partial}$ to be the one-point compactification of $E$. For example, if $E$ is a bounded smooth domain of
$\R^d$, $\xi$ is the killed Brownian motion in $E$ and ${\partial}$ was added as an isolated point, then $\xi$ will not be a Hunt process.
Let \emph{the branching mechanism} $\psi$ be a function on $E \times \mathbb R_+$ given by
\begin{align}
	\psi(x,z)
    = -\beta(x) z + \frac{1}{2}\alpha(x) z^2 + \int_{(0,\infty)} (e^{-zr} -1 + zr) \pi(x,dr),
	\quad x\in E, z\geq 0,
\end{align}
	where $\alpha\geq 0$ and $\beta$ are both in $\mathcal B_b(E)$,  and $\pi$ is a  kernel  from $(E, {\mathcal B}(E))$ to $(\R_+, {\mathcal B}(\R_+))$ satisfying
  \begin{equation}\label{assumption on bm}
\int_0^\infty (r\wedge r^2)\pi(\cdot,dr)\in \mathcal B^+_b(E).
\end{equation}
Note that this assumption implies that, for any fixed $z>0$, $\psi(\cdot,z)$ is bounded on $E$.
We extend $\psi$ to a branching mechanism $\overline \psi$ on $E_\partial$ by defining $\overline \psi(\partial, z)=0$ for all $z\ge 0$.

Let $\mathcal M(E)$  (resp. $\mathcal M(E_\partial)$) denote the space of all finite Borel measures on $E$ (resp. $E_\partial$)
equipped with the topology of weak convergence. Any $\mu\in \mathcal M(E)$ will be identified with its zero extension in
$\mu\in \mathcal M(E_\partial)$. For any $\mu \in \mathcal M(E_\partial)$ and $f\in \mathcal B(E_\partial)$,
we use $\langle f, \mu\rangle$ or $\mu(f)$ to denote the integral of $f$ with respect to $\mu$ whenever the integral is well-defined.
For  $f \in \mathcal B^+_b(E_\partial)$, there is a unique locally bounded non-negative map $(t,x)\mapsto \overline V_tf(x)$ on $\mathbb R_+\times E_\partial$ such that
\begin{equation} \label{eq:M.1-bar}
\overline V_tf(x) + \Pi_x\Big[\int_0^{t} \overline \psi\big(\xi_s, \overline V_{t-s} f(\xi_s)\big) ds\Big] = \Pi_x[f(\xi_t)], \quad t\geq 0, x\in E_\partial.
\end{equation}
	Here, local boundedness of the map $(t,x) \mapsto \overline V_tf(x)$ means that  for any $T>0$,
\[
	\sup_{0\leq t\leq T, x\in E} \overline V_tf(x) < \infty.
\]
Similarly, for  $f \in \mathcal B^+_b(E)$,
there is a unique locally bounded non-negative map $(t,x)\mapsto V_tf(x)$ on $\mathbb R_+\times E$ such that
\begin{equation} \label{eq:M.1}
	V_tf(x) + \Pi_x\Big[\int_0^{t\wedge \zeta} \psi\big(\xi_s, V_{t-s} f(\xi_s)\big) ds\Big] = \Pi_x[f(\xi_t) 1_{\{t< \zeta\}}], \quad t\geq 0, x\in E.
\end{equation}
There exists an $\mathcal M(E_\partial)$-valued Hunt process
$\overline X =\{(\overline X_t)_{t\geq 0}; \mathbb P_\mu, \mu \in \mathcal M(E_\partial)\}$ such that
$$
\mathbb P_\mu[e^{- \overline X_t(f)}]	= e^{- \mu(\overline V_tf)},	\quad t\geq 0, f \in \mathcal B^+_b(E_\partial).
$$
This process $\overline X$ is known as a $(\xi, \overline \psi)$-superprocess.
See \cite[Section 2.3 and Theorem 5.11]{LZ} for more details.
Let $\iota(\mu)$ be the restriction of a measure $\mu\in \mathcal M(E_\partial)$ to $E$ and $X_t=\iota(\overline X_t)$. It follows from the proof of \cite[Theorem 5.12]{LZ} that $X =\{(X_t)_{t\geq 0}; \mathbb P_\mu, \mu \in \mathcal M(E)\}$
is an $\mathcal M(E)$-valued  Markov process such that
$$
\mathbb P_\mu[e^{- X_t(f)}]	= e^{- \mu(V_tf)},	\quad t\geq 0, f \in \mathcal B^+_b(E).
$$
However, since we have taken $E_\partial$ to be the one-point-compatification of $E$, $X$ is in general not a Hunt process and does not have good regularity properties.
Since $X_t(f)=\overline X_t(f)$ for any function $f$ on $E$ and we are only interested in quantities of the form $X_t(f)$, we can work with the Hunt process
$\overline X$ when necessary.
When the initial value is $\delta_x, x\in E$, we write $\mathbb P_x$ for $\mathbb P_{\delta_x}$.
We use $(P^\beta_t)_{t\ge 0}$ to denote the following Feynman-Kac semigroup
$$
P^\beta_tf(x)=
\Pi_x\left(\exp\left(\int^t_0\beta(\xi_s)ds\right)f(\xi_t)I_{\{t<\zeta\}}\right),
\quad x\in E,f \in \mathcal B_b^+(E).
$$
Then it is known (see \cite[Proposition 2.27]{LZ}) that for any $\mu\in\mathcal M(E)$,
\begin{equation} \label{eq:M.2}
	\mathbb P_\mu[X_t(f)] = \mu (P_t^\beta f),
	\quad t\geq 0, f \in \mathcal B_b^+(E).
\end{equation}
$(P_t^\beta)_{t\geq 0}$ is called the mean semigroup  of $X$. For this mean semigroup, we will always assume that
\begin{assumption}\label{asp:H1}
	There exist a constant $\lambda>0$, a positive function $\phi \in \mathcal B_b(E)$ and a probability measure $\nu$ with full support on $E$ such that for any $t\geq 0$, $P_t^\beta \phi = e^{\lambda t} \phi$, $\nu P_t^\beta = e^{\lambda t} \nu$ and $\nu(\phi) =1$.
\end{assumption}

	Denote by $L_1^+(\nu)$ the collection of non-negative Borel functions on $E$ which are integrable with respect to the measure $\nu$.
	Denote by $\mathbf 0$ the null measure on $E$ and $E_\partial$.
	Write $\mathcal M^0(E) = \mathcal M(E)\setminus \{\mathbf 0\}$ and $\mathcal M^0(E_\partial) = \mathcal M(E_\partial)\setminus \{\mathbf 0\}$

	We further assume that
the following assumption holds:
\begin{assumption}\label{asp:H2}
	For all $t>0$, $x\in E$, and $f\in L_1^+(\nu)$, it holds that $P_t^\beta f(x) = e^{\lambda t} \phi(x) \nu(f) (1+ C_{t,x,f})$ for some $C_{t,x,f}\in \mathbb R$, and that $\lim_{t\to\infty}c_t=0$, where $c_t:=\sup_{x\in E, f\in L_1^+(\nu)} |C_{t,x,f}|$.
\end{assumption}
Note that $\lim_{t\to\infty}c_t=0$ implies that there exists $t_0>0$ such that
\[
\sup_{t>t_0}\sup_{x\in E, f\in L_1^+(\nu)} |C_{t,x,f}|
	< \infty.
\]
Without loss of generality,
throughout this paper we will assume $t_0=1$.

For examples satisfying Assumptions 1 and 2, see \cite[Section 1.3]{LRSS} and \cite[Section 1.4]{RSZ}.

Define
 \begin{equation}\label{def-M}
 M_t(\phi):=e^{-\lambda t}\langle\phi, X_t\rangle, \quad t\ge 0.
 \end{equation}
It follows from \eqref{eq:M.2} and Assumption 1 that
$\{M_t(\phi), t\ge 0\}$ is a non-negative c\`{a}dl\`{a}g martingale,
see \eqref{stochatic r for mart} below.
By the martingale convergence theorem, $M_t(\phi)$ has an almost sure
limit as $t\to\infty$. We denote this limit as $M_\infty(\phi)$.
In this paper, we study the rate at which
$M_t(\phi)$ converges to $M_\infty(\phi)$ as $t\to\infty$.

To state our results we need to introduce some notation.
 Define a new kernel
$\pi^\phi(x,dr)$ from $(E, {\mathcal B}(E))$ to $(\R_+, {\mathcal B}(\R_+))$ such that for any non-negative Borel function $f$ on
$\R_+$,
\begin{equation}\label{phi-change}
\int_0^\infty f(r)\pi^\phi(x,dr)=\int_0^\infty f(r\phi(x))\pi(x, dr), \quad x\in E.
\end{equation}
By \eqref{assumption on bm} and the boundedness of
$\phi$, $\pi^\phi$ satisfies
\begin{eqnarray*}
\int_0^\infty (r\wedge r^2)\pi^\phi(x,dr)&=&\phi(x)\Big[\int_0^{1/\phi(x)}r^2\phi(x) \pi(x,dr)+\int_{1/\phi(x)}^\infty r\pi(x,dr)\Big]\\
&\leq& \phi(x)\left(\|\phi\int_0^1r^2\pi(\cdot,dr)\|_\infty+2\|\int_{1\wedge 1/\|\phi\|_\infty}^\infty r\pi(x,dr)\|_\infty\right).
\end{eqnarray*}
Denote by $C:=\|\phi\int_0^1r^2\pi(\cdot,dr)\|_\infty+2\|\int_{1\wedge 1/\|\phi\|_\infty}^\infty r\pi(x,dr)\|_\infty$. Then
\begin{equation}\label{assum: moment finite}
\int_0^\infty (r\wedge r^2) \pi^\phi(x,dr)\leq C\phi(x).
\end{equation}

In \cite{LRS}, we studied the relationship
between $M_\infty(\phi)$ being a non-zero random variable
and the following function $l$:
\begin{equation}\label{m}
l(y):=\int_1^\infty r\ln r \pi^\phi(y, dr),\quad y\in E,
\end{equation}
and established an $L\log L$ criterion
(see Proposition \ref{degeneracy theorem} below)
for a class of superdiffusions
with $\alpha=0$.  It is easy to check that this criterion still holds for
the superprocesses described above.

\begin{prop}\label{degeneracy theorem}\cite[Theorem 1.1]{LRS}
Suppose that Assumptions \ref{asp:H1}-\ref{asp:H2} hold and $\mu\in \mathcal M^0(E)$.
Then
$\mathbb P_\mu (M_\infty(\phi))=\langle\phi,\mu\rangle$ if and only if the following $L\log L$ condition holds:
\begin{equation}\label{LlogL}
\int_El(y)\nu(dy)<\infty.
\end{equation}
Moreover, if \eqref{LlogL} holds, then for any $\mu\in \mathcal M_f^0(E)$,
$$
\{M_\infty(\phi)>0\}=\{X_t>0,\ \forall t>0\}, \quad \mathbb P_\mu\mbox{-a.s.}
$$
Otherwise,  $M_\infty(\phi)=0$, $\mathbb P_\mu$-a.s. for any $\mu\in \mathcal M^0(E)$.
\end{prop}

Throughout this paper, we assume that \eqref{LlogL} holds. Thus $M_t(\phi)$ converges
to $M_\infty(\phi)$ $\mathbb P_\mu$-almost surely and in $L^1(\mathbb P_\mu)$  for any $\mu\in \mathcal M^0(E)$.

Since $M_s(\phi)$ is right continuous, for any $a^*>0$, we can define
\begin{equation}\label{def: A(a)}
A_t(a^*)=\int_0^t e^{\frac{\lambda s}{a^*}}\big(M_\infty(\phi)-M_s(\phi)\big)ds,\qquad t\in[0,\infty).
\end{equation}
Note that
$$
A_t(a^*)-A_1(a^*)=\int_1^t e^{\frac{\lambda s}{a^*}}\big(M_\infty(\phi)-M_s(\phi)\big)ds,\quad t\geq 1.
$$
The convergence of $A_t(a^*)-A_1(a^*)$ as $t\to\infty$ is related to the
rate at which $M_\infty(\phi)-M_t(\phi)$ converges to $0$ as $t\to\infty$.
Our first result is the following criterion for the
$L^p$ convergence rate of $M_\infty(\phi)-M_t(\phi)$ to $0$ as $t\to\infty$.
We use the usual notation $\|\cdot\|_p$ to denote the $L^p$ norm with $p\ge 1$.

\begin{thm}\label{theorem:Lpconv}
  Assume that Assumptions \ref{asp:H1}-\ref{asp:H2}  and \eqref{LlogL} hold. Let $1<a<p\leq 2$ and $\frac1a+\frac1{a^*}=1$.
\begin{itemize}
\item[(1)] If
\begin{equation}\label{finite of p moment}
\int_E\nu(dy)\int_1^\infty r^p \pi^\phi(y, dr)<\infty,
\end{equation}
then for any $\mu\in\mathcal M^0(E)$, $(A_t(a^*)-A_1(a^*))$ converges in $L^p(\mathbb P_\mu)$ and $\mathbb P_\mu$-almost surely as $t\to\infty$.
\item[(2)] If for some $\mu\in\mathcal M^0(E)$, $(A_t(a^*)-A_1(a^*))$ converges in $L^p(\mathbb P_\mu)$ as $t\to \infty$,
then it must converge $\mathbb P_\mu$-almost surely and \eqref{finite of p moment} holds.
\item[(3)] If \eqref{finite of p moment} holds, $\left\|M_\infty(\phi)-M_t(\phi)\right\|_p=o(e^{-\frac{\lambda t}{a^*}})$ as $t\to\infty$.
\item[(4)]If $\left\|M_\infty(\phi)-M_t(\phi)\right\|_p=o(1)$ as $t\to\infty$, then \eqref{finite of p moment} holds.
\end{itemize}
\end{thm}

For a Galton-Watson process $Z$, it is proved in \cite[Proposition 1.3]{Liu01} that if $E(Z_1^p)<\infty$ for some $p>1$, then there exists some $c>0$ such that
$$
\|W_n-W_\infty\|_p\le \left\{\begin{array}{ll}c m^{-\frac{1}{q}n}, \quad &\mbox{if } p\in(1,2], \\
cm^{-\frac{1}{2}n},\quad &\mbox{if } p>2.\end{array}\right.
$$
The above results imply that $\|W_n-W_\infty\|_p=o(m^{-\frac{n}{a^*}})$ for $1<a<p\le 2$, which corresponds to our Theorem \ref{theorem:Lpconv}(3), and  $\|W_n-W_\infty\|_p=o(\rho^{-n})$ for any $\rho<m^{1/2}$.

If $\int^\infty_1r^2\pi(x, dr)$ is bounded (which implies that \eqref{finite of p moment} holds for $p=2$),
by  the central limit theorem
 (see \cite[Theorem 1.4]{RSZ2}),
$e^{-\lambda t/2}(M_t(\phi)-M_\infty(\phi))$ converges to $ Z\sqrt{M_\infty(\phi)}$  with $Z$ being a normal random variable with mean zero and independent of $M_\infty(\phi)$. In \cite{RSZ2}, the  mean semigroup $P^\beta_t$ is assumed to be  symmetric  with respect to some measure $m$, and the assumptions on $(P^\beta_t)_{t\geq 0}$ are slightly different, but the central limit theorem also holds in the nonsymmetric  case, see \cite{RSZ3} for the corresponding results for branching Markov processes.
A question related to Theorem \ref{theorem:Lpconv} is whether the results still holds
for $a=p<2$.
The following theorem gives necessary and sufficient conditions for
the almost sure convergence for the case of
$a=p<2$.
\begin{thm}\label{thm: p moment convergence rate}
Suppose that Assumptions \ref{asp:H1}-\ref{asp:H2} and \eqref{LlogL} hold.
Let $1<p<2$, $1/p+1/q=1$.
\begin{itemize}
\item[(1)] If \eqref{finite of p moment} holds, then for any $\mu\in \mathcal{M}^0(E)$,
as $t\to\infty$, $A_t(q)$ converges $\mathbb P_\mu$-a.s. and
\[
M_t(\phi)-M_\infty(\phi)=
o(e^{-\frac{\lambda t}{q}}),
\quad\mathbb P_\mu\mbox{-a.s.}
\]
\item[(2)] Suppose there exist $B>0$ and $T_0>0$ such that
 \begin{equation}\label{assm: unif upp}
\sup_{x\in E}\dfrac{1}{\phi(x)}\int_t^\infty \pi^\phi(x,dr)\leq B \int_E\nu(dy)\int_t^\infty \pi^\phi(y,dr),\qquad   t>T_0.
 \end{equation}
If
\begin{equation}\label{m2}
\int_E\nu(dy)\int_1^\infty r^p \pi^\phi(y, dr)=\infty,
\end{equation}
then  for any $\mu\in \mathcal{M}^0(E)$,
$ M_t(\phi)-M_\infty(\phi)=
o(e^{-\frac{\lambda t}{q}})$
$\mathbb P_\mu\mbox{-a.s.}$ does not hold as $t\to\infty$.
\end{itemize}
\end{thm}

\begin{thm}\label{thm: log}
Assume  that Assumptions \ref{asp:H1}-\ref{asp:H2}
and \eqref{LlogL} hold.
\begin{itemize}
\item[(1)]For any $\gamma> 0$,
\begin{equation}\label{assum: log moment}
\int_E\nu(dx)\int_1^\infty r(\ln r)^{\gamma+1}\pi^\phi(x,dr)<\infty
\end{equation}
implies that, for any $\mu\in \mathcal{M}^0(E)$,
$$
\int_0^t s^{\gamma-1}(M_\infty(\phi)-M_s(\phi))ds
$$
 converges $\mathbb P_\mu$-almost surely as $t\to\infty$, and
\[
M_\infty(\phi)-M_t(\phi)=o(t^{-\gamma}),\qquad \mathbb P_\mu\mbox{-a.s.}
\]
If \eqref{assum: log moment} holds with $\gamma\geq 1$, then
$\int_0^\infty (M_\infty(\phi)-M_t(\phi))dt$ exists $\mathbb P_\mu$-almost surely
for any $\mu\in \mathcal{M}^0(E)$.
\item[(2)]Suppose that there exist $b>0, T_1>0$ and
a Borel set $F\subset E$  with $\nu(F)>0$ such that
 \begin{equation}\label{assum: lower}
 \inf_{x\in F}\dfrac{1}{\phi(x)}\int_t^\infty r\pi^\phi(x,dr)\geq b\int_E\nu(dx)\int_t^\infty r\pi^\phi(x,dr),
 \qquad t>T_1.
 \end{equation}
 If there is $\gamma\in (0,\infty)$ such that
\begin{equation}\label{assum: log moment infty}
\int_E\nu(dx)\int_1^\infty r(\ln r)^{\gamma+1}\pi^\phi(x,dr)=\infty,
\end{equation}
then
for any $\mu\in \mathcal{M}^0(E)$,
$\int_0^t s^{\gamma-1}(M_\infty(\phi)-M_s(\phi))ds$
does not converge $\mathbb P_\mu\mbox{-a.s.}$ as $t\to\infty$.
If, as $t\to\infty$,
\begin{equation}\label{inf of log}
\int_E\nu(dx)\int_t^\infty r(\ln r-\ln t)\pi^\phi(x,dr)=
o\left((\ln t)^{-\gamma}\right)
\end{equation}
does not
hold, then $M_\infty(\phi)-M_t(\phi)=o(t^{-\gamma}),\,\mathbb P_\mu\mbox{-a.s.}$ does not hold as well.
\end{itemize}
\end{thm}
It was noted   in \cite[Theorem II.4.1, p.36]{AH}
that \eqref{assum: log moment} implies \eqref{inf of log}, and \eqref{inf of log}
 implies  that
 $$\int_E\nu(dx)\int_1^\infty r(\ln r)^{\gamma+1-\varepsilon}\pi^\phi(x,dr)<\infty\quad \mbox{ for all } 0<\varepsilon\leq \gamma.$$
This says that \eqref{assum: log moment} is slightly stronger than \eqref{inf of log}.

We make a few remarks about \eqref{assm: unif upp} and \eqref{assum: lower}.
Note that by definition,
$$
\int^\infty_t\pi^\phi(x, dr)=\int^\infty_{t/\phi(x)}\pi(x, dr)=\pi\left(x,(t/\phi(x), \infty)\right).
$$
If $\pi(x, dr)=\gamma(x)r^{-1-\alpha}dr$ with $\alpha\in (1, 2)$ and $\gamma$ a bounded non-negative Borel function, then
$$
\int^\infty_t\pi^\phi(x, dr)=\frac1\alpha\gamma(x)t^{-\alpha}\phi(x)^\alpha.
$$
Hence
$$
\frac1{\phi(x)}\int^\infty_t\pi^\phi(x, dr)=\frac1\alpha\gamma(x)t^{-\alpha}\phi(x)^{\alpha-1}
$$
and
$$
\int_E\nu(dx)\int^\infty_t\pi^\phi(x, dr)=\frac{t^{-\alpha}}{\alpha} \int_E\nu(dx)\gamma(x)\phi(x)^\alpha.
$$
Since $\gamma$ and $\phi$ are bounded, \eqref{assm: unif upp} is satisfied. Similarly, by definition,
$$
\int^\infty_tr\pi^\phi(x, dr)
=\phi(x)\int^\infty_{t/\phi(x)}r\pi(x, dr).
$$
Similarly, we have
$$
\frac1{\phi(x)}\int^\infty_tr\pi^\phi(x, dr)=\frac1{\alpha-1}\gamma(x)\phi(x)^{\alpha-1}t^{1-\alpha}
$$
and
$$
\int_E\nu(dx)\int_t^\infty r\pi^\phi(x,dr)=\frac{t^{1-\alpha}}{\alpha-1}\int_E\gamma(x)\phi(x)^{\alpha}\nu(dx).
$$
 Thus \eqref{assum: lower} is satisfied if $\nu(\{x: \gamma(x)>0\})>0$. It is easy to generalize the remarks above on \eqref{assm: unif upp} and \eqref{assum: lower}
to the case when $\pi(x, dr)=\gamma(x)r^{-1-\alpha}s(r)dr$ with $\alpha\in (1, 2)$, $\gamma$ a bounded non-negative Borel function,
and  $s$ a local
bounded non-negative Borel function $(0, \infty)$ which is slowly varying at $\infty$.

The remainder of this article is organized as follows.  In Section \ref{SIR}, we present a stochastic integral
representation of superprocesses which will be used in later sections.
In Section \ref{spine}, we introduce a spine decomposition of superprocesses
which is used in the proof of Lemma \ref{le: upper}.
The main results are proved in Section \ref{proof-of-results}.
Lemma \ref{le: upper} plays a key role in  the proof of Theorem \ref{thm: p moment convergence rate}.

In this paper, we use the convention that  an expression of the type $a\lesssim b$ means that there exists a positive constant $N$ which is independent of $a$ and $b$ such that $a\leq Nb$.
Moreover, if $a\lesssim b$ and $b\lesssim a$, we shall write $a \asymp b$.

\section{Superprocesses}
\subsection{Stochastic Integral Representation of Superprocesses}\label{SIR}
Without loss of generality, we assume that our process
$\overline X$ is the coordinate process on
\[
  \mathbb D:=\{ w= (w_t)_{t\geq 0}: w \text{ is an
   $\mathcal M(E_\partial)$-valued c\`{a}dl\`{a}g function on $[0,\infty)$.}
 \}.
\]
We assume that
$(\mathcal{F}_\infty, (\mathcal{F}_t)_{t\ge 0})$ is
the natural filtration on $\mathbb D$, completed as usual with the $\mathcal{F}_\infty$-measurable and $\mathbb P_\mu$-negligible sets
for every $\mu\in\mathcal{M}(E_\partial)$.
Let $\mathbb W^+_0$ be the family of
$\mathcal M(E_\partial)$-valued c\`{a}dl\`{a}g
functions on $(0, \infty)$ with $\mathbf 0 $ as a trap and with
$\lim_{t\downarrow0}w_t= \mathbf 0$. $\mathbb W^+_0$ can be regarded as a subset of $  \mathbb D$.

Throughout this paper  assume that $\mathbb P_{x}(X_t(1)=0)>0$ for any $x\in E$ and $t>0$, which implies that  there exists a unique family of
$\sigma$-finite measures $\{\mathbb N_x; x\in E\}$ on
$\mathbb W^+_0$
such that
for any $\mu\in \mathcal M(E)$, if ${\mathcal N}(dw)$ is a Poisson random measure on
$\mathbb W^+_0$ with intensity measure
$$
\mathbb N_\mu(dw):=\int_E \mathbb N_x(dw)\mu(dx),
$$
then the process defined by
$$
\widetilde X_0=\mu, \quad \widetilde X_t=\int_{\mathbb W^+_0}w_t{\mathcal N}(dw), \quad t>0
$$
is a realization of the superprocess
$\overline X=\{(\overline X_t)_{t\geq 0}; \mathbb P_\mu, \mu \in \mathcal M(E)\}$.
Furthermore,
$ \mathbb N_x(\langle f,w_t\rangle)=
\mathbb P_{x}\langle f, X_t\rangle$ for any $f\in \mathcal B^+(E)$
(see \cite[Theorem 8.22]{LZ} and \cite[Section 2.2]{RSZ4}).
$\{\mathbb N_x; x\in E\}$ can be regarded as measures on $\mathbb D$ carried by
$\mathbb W^+_0$.

Let us recall the stochastic integral representation of
superprocesses, for more details see \cite{F2} or \cite[Chapter 7]{LZ}.
Let $({\bf A}, \mathfrak{D}({\bf A}))$  be the weak infinitesimal generator of $\xi$ as defined in \cite[Section 4]{F1}.
For any $f\in \mathfrak{D}({\bf A})$,
$$
\frac{P_tf(x)-f(x)}{t}\to {\bf A}f(x), \quad\mbox{bounded and pointwisely as }t\to 0.
$$
We will use the standard notation
$\triangle \overline X_s=\overline X_s-\overline X_{s-}$ for the jump of $\overline X$ at time $s$.
Let $C_0^2(\mathbb{R})$ denote the set of all twice continuously differentiable functions on $\R$ vanishing at infinity.
It is known (cf. \cite[Theorem 1.5]{F1})
that the superprocess $\overline X$ is a
solution to the following martingale problem:
for any $\varphi \in \mathfrak{D}({\bf A})$
and $h\in C_0^2(\R)$,
\begin{equation}\label{eq:martingale problem}
\begin{split}
&\displaystyle h(\langle\varphi, \overline X_t\rangle)-h(\langle\varphi,\mu\rangle)
-\int_0^th'(\langle\varphi, \overline X_s \rangle)\langle
({\bf A}+\beta)\varphi, \overline X_s\rangle ds
-\frac{1}{2}\int_0^th''(
\langle\varphi, \overline X_s \rangle)\langle\alpha\varphi^2, \overline X_s\rangle ds\\
&\displaystyle-\int_0^t\int_E\int_{(0, \infty)}\big(h(\langle \varphi,
\overline X_s\rangle+r\varphi(x))-h(\langle \varphi, \overline X_s\rangle)-h'(\langle
\varphi, \overline X_s\rangle)r\varphi(x)\big)\pi(x, dr) \overline X_{s-}(dx)ds
\end{split}
\end{equation}
is a $\mathbb P_\mu$-martingale
 for any $\mu\in \mathcal M^0(E_\partial)$.

By  \cite[Proposition 2.1]{F2} (also see \cite[Theorem 7.13]{LZ}),
for any $\varphi\in \mathfrak{D}({\bf A})$ and
$\mu\in \mathcal M^0(E_\partial)$,
\begin{equation}\label{Martingale-expres1}
\langle \varphi, \overline X_t\rangle = \langle\varphi, \mu\rangle+S^J_t(\varphi)+ S^C_t(\varphi) + \int_0^t\langle({\bf A}+\beta)\varphi, \overline X_s\rangle ds,
\end{equation}
where
$ S^C_t(\varphi)$ is a continuous $\mathbb P_\mu$-local martingale and
$S^J_t(\varphi)$ is a $\mathbb P_\mu$-pure jump martingale.
The quadratic variation process of the continuous local martingale $S^C_t(\varphi)$ is given by
\begin{equation}\label{quadratic variation}
\langle S^C(\varphi)\rangle_t=\int_0^t\langle \alpha\varphi^2,
\overline X_s\rangle ds.
\end{equation}

Next, we characterize the pure jump martingale $(S^J_t(\varphi), t\geq 0)$.
Let $J$ denote the set of all jump times of $\overline X$ and
$\delta$ denote the Dirac measure.
From the last part of \eqref{eq:martingale problem}, we see that
the only possible jumps of $\overline X$ are point measures of the form
$r\delta_x$ with $r>0$ and $x\in E_\partial$, see \cite[Section 2.3]{LM}.
Thus the predictable compensator of the random measure (for the definition of the predictable compensator of a random measure, see, for instance \cite[p.107]{D})
$$
N:=\sum_{s\in  J}\delta_{(s, \triangle \overline X_s)}
$$
is a random measure $\widehat{N}$ on $\R_+\times \mathcal{M}(E_\partial)$
such that for any nonnegative predictable function $F$ on
$\R_+\times \Omega\times \mathcal{M}(E_\partial)$,
\begin{eqnarray}\label{poisson point definition}
\int_0^\infty\int_{\mathcal M(E_\partial)} F(s,\omega, \upsilon)\widehat{N}(ds, d\upsilon)
=\int_0^\infty ds\int_E \overline X_{s-}(d x) \int_0^\infty
F(s,\omega, r\delta_x)\pi(x,dr),
\end{eqnarray}
where $\pi(x,dr)$ is the kernel of the branching mechanism
$\psi$.  Therefore we have
\begin{eqnarray}\label{expectation identity}
\mathbb P_\mu\left[\sum_{s\in J}F(s,\omega,\triangle \overline X_s)\right]
=\mathbb P_\mu\int_0^\infty ds\int_E \overline X_{s-}(dx)
\int_0^\infty  F(s,\omega, r\delta_x)\pi(x,dr).
\end{eqnarray}
 See \cite[p.111]{D}.

Let $F$ be a Borel function on
$\mathbb{R}_+\times \mathcal{M}(E_\partial)$ satisfying
$$
\mathbb P_\mu\bigg[\bigg(\sum_{s\in[0,t],s\in J}F(s, \Delta
\overline X_s)^2\bigg)^{1/2}\bigg]<\infty, \qquad \mbox{ for all } \mu\in
{\mathcal M}(E_\partial).
$$
Then the stochastic integral of $F$ with respect to the compensated
random measure $N-\widehat{N}$
$$
\int_0^t\int_{\mathcal{M}(E_\partial)}F(s, \upsilon)(N-\widehat{N})(ds, d\upsilon)
$$
can be defined (cf. \cite{LM} and the references therein) as the
unique purely discontinuous martingale (vanishing at time $0$) whose
jumps are indistinguishable from
$1_J(s)F(s,\Delta \overline X_s)$.

Suppose that $g$ is a Borel function on $\R_+\times E$.  Define
\begin{eqnarray}
F_g(s,\upsilon):=\int_Eg(s,x)\upsilon(dx),\quad \upsilon\in \mathcal{M}(E_\partial),
\end{eqnarray}
whenever the integral above makes sense.
When $g$ is  a bounded Borel function on $\mathbb R_+\times E$, for any $\mu\in \mathcal{M}(E)$,
\begin{equation*}
\begin{array}{rl}
&\mathbb P_\mu\left[\left(\sum_{s\in[0,t],s\in J} F_g(s, \Delta \overline X_s)^2\right)^{1/2}\right]=\mathbb P_\mu\left[\left(\sum_{s\in[0,t],s\in J}
\left(\int_E g(s,x)(\Delta \overline X_s)(dx)\right)^2\right)^{1/2}\right]\\
\le&\displaystyle\|g\|_\infty \mathbb P_\mu\left[\left(
\sum_{s\in[0,t],s\in J}\langle 1, \Delta \overline X_s\rangle^2 I_{\{\langle
1, \Delta \overline X_s\rangle\le 1\}}+\sum_{s\in[0,t],s\in J} \langle 1,
\Delta \overline X_s\rangle^2 1_{\{\langle 1, \Delta \overline X_s\rangle> 1\}}
\right)^{1/2}\right]\\
\le&\displaystyle\|g\|_\infty \mathbb P_\mu\left[\left(
\sum_{s\in[0,t],s\in J}\langle 1, \Delta \overline X_s\rangle^2 1_{\{
\langle 1, \Delta \overline X_s\rangle\le 1\}}\right)^{1/2}\right]\\
&+\displaystyle\|g\|_\infty \mathbb P_\mu\left[\left(
\sum_{s\in[0,t],s\in J}\langle 1, \Delta \overline X_s\rangle^2 1_{\{ \langle
1, \Delta \overline X_s\rangle> 1\}}\right)^{1/2}\right].
\end{array}
\end{equation*}
Using the first two displays on \cite[p.203]{LM}, we get
 \begin{eqnarray}\label{int-finit}
\mathbb P_\mu\left[\left(\sum_{s\in[0,t],s\in J}F_g(s, \Delta \overline X_s)^2\right)^{1/2}\right]<\infty.
\end{eqnarray}
Therefore, if $g$ is bounded on $\R_+\times E$, then the integral $\int_0^t\int_{\mathcal{M}(E_\partial)}F_g(s,\upsilon)(N-\widehat{N})(ds,d\upsilon)$ is well defined and is a martingale.  Define the martingale measure $S^J(ds,dx)$ by
\begin{eqnarray}\label{martingale definition of SJ}
\int_0^t\int_Eg(s, x)S^J(ds,dx):=\int_0^t\int_{\mathcal{M}(E_\partial)}F_g(s,\upsilon)(N-\widehat{N})(ds,d\upsilon).
\end{eqnarray}
Thus the pure jump martingale $S^J_t(\varphi)$ in \eqref{Martingale-expres1} can be written as
\[
S^J_t(\varphi)=\int_0^t\int_E\varphi(x)S^J(ds,dx).
\]

A martingale measure $S^C(ds,dx)$
can be defined (see \cite{F2} or \cite{P} for the precise definition) so
that the continuous martingale in \eqref{Martingale-expres1} can be expressed as
\[
 S_t^C(\varphi)=\int_0^t\int_{E}\varphi(x)S^C(ds,dx).
\]
Summing up these two martingale measures, we get a martingale measure
\begin{equation}\label{M=J+C}
M(ds,dx)=S^J(ds,dx)+S^C(ds,dx).
\end{equation}

Using this, \cite[Proposition 2.14]{F2} and applying a limit argument,
one can show that for any bounded Borel function $g$ on $E$,
\begin{equation}\label{dual representation}
\langle g, X_t\rangle=\langle P^{\beta}_tg,\mu\rangle+\int_0^t\int_E
P^{\beta}_{t-s}g(x)S^J(ds,dx)+\int_0^t\int_E
P^{\beta}_{t-s}g(x)S^C(ds,dx).
\end{equation}
In particular, taking $g=\phi$ in \eqref{dual representation}, where
$\phi$ is the positive eigenfunction of $P^{\beta}_t$ given in Assumption \ref{asp:H1},
 we get the expression for the martingale $(M_t(\phi))_{t\geq 0}$:
\begin{eqnarray}\label{stochatic r for mart}
M_t(\phi)=\langle \phi,\mu\rangle+\int_0^t e^{-\lambda s}\int_E \phi(x)S^J(ds, dx)+\int_0^t e^{-\lambda s} \int_E\phi(x)S^C(ds, dx).
\end{eqnarray}
Therefore the limit $M_\infty(\phi)$ of $ M_t(\phi)$ can be written as
\begin{equation}\label{martingale representation of limit}
M_\infty(\phi)=
\langle \phi,\mu\rangle+\int_0^{\infty} e^{-\lambda s}\int_E \phi(x)S^J(ds,dx)+\int_0^{\infty} e^{-\lambda s}\int_E \phi(x)S^C(ds, dx).
\end{equation}

 For the jump part above, we always handle the `small jumps' and the `large jumps' separately. Now let us give the precise definitions of `small jumps' and `large jumps'. Given $\rho\in(0,\infty]$, a jump at time $s$ is called `small' if $0<\Delta \overline X_s(\phi)<e^{\frac{\lambda}{\rho} s}$,  and
`large' if $\triangle \overline X_s(\phi)\ge e^{\frac{\lambda}{\rho} s}$,
where $\triangle \overline X_s(\phi)=r\phi(x)$ when $\triangle \overline X_s=r\delta_x$ with $r>0$ and $x\in E$.  Define
$$
N^{(1,\rho)}:=\sum_{0<\triangle \overline X_s(\phi)<e^{\frac{\lambda}{\rho}s}}\delta_{(s, \triangle \overline X_s)}, \quad
N^{(2,\rho)}:=\sum_{\triangle \overline X_s(\phi)\ge e^{\frac{\lambda}{\rho} s}}\delta_{(s, \triangle \overline X_s)},
$$
and denote the compensators of $N^{(1,\rho)}$  and $N^{(2,\rho)}$ by
$\widehat N^{(1,\rho)}$ and $\widehat N^{(1,\rho)}$ respectively.
Then for any non-negative Borel function $F$ on $\mathbb{R}_+\times \mathcal{M}(E_\partial )$,
$$
\int_0^\infty \int_{\mathcal{M}(E_\partial)} F(s, \upsilon)\widehat N^{(1,\rho)}(ds, d\upsilon)
=\int_0^\infty ds\int_E\overline X_{s-}(dx)\int_0^{e^{\frac{\lambda}{\rho}s}} F(s, r\phi(x)^{-1}\delta_x)\pi^{\phi}(x,dr)
$$
and
$$
\int_0^\infty \int_{\mathcal{M}(E_\partial)} F(s,\upsilon)\widehat N^{(2,\rho)}(ds, d\upsilon)
=\int_0^\infty ds\int_E\overline X_{s-}(dx)\int_{e^{\frac{\lambda}{\rho}s}}^\infty F(s, r\phi(x)^{-1}\delta_x)\pi^{\phi}(x, dr),
$$
where $\pi^\phi$ was defined in \eqref{phi-change}.
Let $J^{(1,\rho)}$ denote the set of jump times of $ N^{(1,\rho)}$, and let $J^{(2,\rho)}$ denote the set of jump times of $ N^{(2,\rho)}$.
Then
\begin{equation*}\label{identity1a}
\int_0^\infty \int_{\mathcal{M}(E_\partial )} F(s,\upsilon) N^{(1,\rho)}(ds, d\upsilon)=\sum_{s\in J^{(1,\rho)}}F(s,\triangle \overline X_s),
\end{equation*}
\begin{equation*}\label{identity1(2)a}
\int_0^\infty \int_{\mathcal{M}(E_\partial )} F(s,\upsilon) N^{(2,\rho)}(ds,d\upsilon)=\sum_{s\in J^{(2,\rho)}}F(s,\triangle \overline X_s).
\end{equation*}
Similar to the way we
 constructed $S^J(ds,dx)$ from $N(ds,d\upsilon)$,
we can construct two martingale measures $S^{(1,\rho)}(ds,dx)$ and $S^{(2,\rho)}(ds, dx)$
respectively from $N^{(1,\rho)}(ds,d\upsilon)$ and $N^{(2,\rho)}(ds, d\upsilon)$.  Then for any bounded Borel function $g$ on $\R_+\times E$, we can obtain the following martingales, for $t>0$,
\begin{equation}\label{martingale-int2}
S^{(1,\rho)}_t(g)=\int^t_0\int_Eg(s,x) S^{(1,\rho)}(ds,dx)=
\int^t_0\int_{{\mathcal M}(E_\partial)}F_g(s,\upsilon)(N^{(1,\rho)}-\widehat
N^{(1,\rho)})(ds, d\upsilon)
\end{equation}
and
\begin{equation}\label{martingale-int2(2)}
S^{(2,\rho)}_t(g)=\int^t_0\int_Eg(s,x) S^{(2,\rho)}(ds, dx)
=\int^t_0\int_{{\mathcal M}(E_\partial)}F_g(s,\upsilon)(N^{(2,\rho)}-\widehat N^{(2,\rho)})(ds, d\upsilon),
\end{equation}
where $F_g(s,\upsilon)=\int_E g(s,x)\upsilon(dx)$.

\subsection{Spine Decomposition of Superprocesses}\label{spine}
 Recall that $\{(\xi_t)_{t\geq 0}; \Pi_x,x\in E\}$ is the spatial motion. Let  $(\mathcal F_t^{\xi})_{t\geq 0}$ be the natural filtration of $(\xi_t)_{t\geq 0}$.  For each $x\in E$, let $\widetilde \Pi_{x}$ be the probability measure defined by
\begin{align}
	\dfrac{d\widetilde{\Pi}_x|_{\mathcal F^{\xi}_t}}{d\Pi_x|_{\mathcal F^{\xi}_t}}=
 \frac{e^{\int_0^t \beta(\xi_s)ds}\phi(\xi_t) \mathbf 1_{\{t<\zeta\}}}{e^{\lambda t}\phi(x)},
	\quad t\geq 0.
\end{align}
It can be verified (see \cite{KS1} for example) that the process
$\{(\xi_t)_{t\geq 0}; \widetilde\Pi_x, x\in E\}$
is a time homogeneous Markov process.
For any $\mu \in \mathcal M(E)$, define $(\phi\mu)(dx):=\phi(x)\mu(dx)$. For any $\mu \in \mathcal M^0(E)$, we define
\[
\Pi_{\mu}(\cdot):= \mu(E)^{-1}\int_{E} \Pi_x(\cdot)\mu(dx)\qquad    \mbox{and} \qquad \widetilde\Pi_{\mu}(\cdot):= \mu(E)^{-1} \int_E\widetilde\Pi_x(\cdot)\mu(dx).
\]
For any $\mu \in \mathcal M^0(E)$, we define the probability measure $\widetilde{\mathbb P}_\mu$ by
\begin{eqnarray}\label{doob trans}
\dfrac{d\widetilde{\mathbb P}_\mu|_{\mathcal F_t}}{d\mathbb P_\mu|_{\mathcal F_t}}=\dfrac{M_t(\phi)}{\mu(\phi)},\quad  t\geq 0.
\end{eqnarray}
It is known (see, for instance, \cite{LRS}) that for any $t>0$,
  \begin{align}\label{spine-decom2}
	\Big( X_t; \widetilde{\mathbb P}_\mu \Big)
	\overset{d}{=}
\Big( X_t + \sum_{\sigma \in \mathcal D^C \cap [0, t]} X^{C,\sigma}_{t-\sigma}+\sum_{\sigma \in \mathcal D^J \cap [0, t]} X^{J,\sigma}_{t-\sigma} ;\, \mathbb Q_\mu\Big).
\end{align}
The right-hand side is constructed as follows.
\begin{itemize}
\item[(i)]
   $\{(\xi_t)_{t\geq 0}; \mathbb Q_\mu\}$ is a Markov process,  called the spine process, with
$$
\{(\xi_t)_{t\ge 0}; \mathbb Q_\mu\}\overset{d}{=}\{(\xi_t)_{t\ge0}; \widetilde \Pi_{\phi\mu}\};
$$
\item[(ii)]
 Conditioned on $(\xi_t)_{t\geq 0}$, \emph{the continuum immigration}
  $\{ (X^{C, \sigma})_{\sigma \in \mathcal D^\mathrm C};
 \mathbb Q_\mu(\cdot |(\xi_t)_{t\geq 0})\}$ is a $\mathbb D$-valued point process such that
\begin{equation}\label{compo con}
	\mathrm n(ds,dw) := \sum_{\sigma\in \mathcal D^C} \delta_{(\sigma, X^{C, \sigma})}(ds,dw)
\end{equation}
is a Poisson random measure on $\mathbb R_+\times \mathbb D$ with intensity
\[
\mathbf n(ds,dw):= \alpha(\xi_s) ds \cdot \mathbb N_{\xi_s}(dw);
\]
\item[(iii)]Conditioned on $(\xi_t)_{t\geq 0}$, \emph{the discrete immigration} $\{(X^{J, \sigma})_{\sigma\in \mathcal D^J}; \mathbb Q_\mu(\cdot |(\xi_t)_{t\geq 0})\}$ is a $\mathbb D$-valued point process such that $\mathrm m({\mathrm d}s,{\mathrm d}w) := \sum_{\sigma\in \mathcal D^J} \delta_{(\sigma, X^{J, \sigma})}(ds,dw)$ is a Poisson random measure on $\mathbb R_+ \times \mathbb D$ with intensity
\begin{align}
 \mathbf m(ds,dw):= ds \cdot \int_{(0,\infty)} y \mathbb P_{y\delta_{\xi_s}}(X\in dw) \pi(\xi_s,dy);
\end{align}
\item[(vi)] 	Given $(\xi_t)_{t\geq 0}$,  $(X^{C,\sigma})_{\sigma \in \mathcal D^C}$ and  $(X^{J,\sigma})_{\sigma\in \mathcal D^J}$ are independent.
\item[(v)]$\{(X_t)_{t\geq 0}; \mathbb Q_\mu\}$ is a copy of the superprocess $\{(X_t)_{t\geq 0}; \mathbb P_\mu\}$, and is independent of  $(\xi_t)_{t\geq 0}$,  $(X^{C, \sigma})_{\sigma \in \mathcal D^C}$ and  $(X^{J, \sigma})_{\sigma\in \mathcal D^J}$.
\end{itemize}
$\{(\xi_t)_{t\geq 0}, (X^{C, \sigma})_{\sigma\in \mathcal D^C}, (X^{J, \sigma})_{\sigma \in \mathcal D^J}, (X_t)_{t\geq 0}; \mathbb Q_\mu\}$ is
called a spine decomposition of $\{(X_t)_{t\geq 0}; \widetilde{\mathbb P}_\mu\}$.

Put
 $$
 Z_t^C:=\sum_{\sigma \in \mathcal D^C \cap [0, t]} X^{C,\sigma}_{t-\sigma}\quad\mbox{ and }
 \quad Z_t^J:=\sum_{\sigma \in \mathcal D^J \cap [0, t]} X^{J,\sigma}_{t-\sigma},\qquad t>0.
 $$
 Then the spine representation \eqref{spine-decom2} of $X$ can be simplified as for any $t\ge 0$,
\begin{align}\label{spine-decom}
	\Big( X_t; \widetilde{\mathbb P}_\mu \Big)
	\overset{d}{=}
\Big( X_t +Z_t^C +Z_t^J ;\, \mathbb Q_\mu\Big).
\end{align}

\section{Proofs of Main Results}\label{proof-of-results}

\subsection{Some Lemmas}
Recall the definition \eqref{def: A(a)}, i.e.,
$$
A_t(q)=\int_0^te^{\frac{\lambda s}{q}}(M_\infty(\phi)-M_s(\phi))ds, \quad t\in [0,\infty).
$$
Let $A(q)$ denote the almost sure limit of $A_t(q)$ as $t\to\infty$  whenever it exists.  For any $p>0$, $g(s,x):=e^{\frac{-\lambda s}{p}}\phi(x)$ is bounded on $\mathbb R_+\times E$, and thus we can define a martingale $(\widetilde{A}_t(p))_{t\ge 0}$ by
\begin{equation}\label{def: ano mart}
\widetilde A_t(p)=\int_0^t e^{\frac{-\lambda s}{p}}\int_E\phi(x)M(ds,dx),\qquad t\geq 0.
\end{equation}
When the almost sure limit of this martingale exists as $t\to\infty$, we denote
the limit by $\widetilde A(p)$ and write it in the integral form
$\widetilde A(p):=\int_0^\infty e^{\frac{-\lambda s}{p}}\int_E\phi(x)M(ds,dx)$.

\begin{lemma}\label{le: equi A A'} Assume that \eqref{LlogL} holds.
Suppose $p\in (1,2]$, $1/p+1/q=1$,
$r> 1$,
and $\mu\in \mathcal{M}^0(E)$.

$(1)$
$\widetilde A_t(p)$  converges $\mathbb{P}_\mu$-almost
surely as $t\to\infty$ if and only if $A_t(q)$ converges $\mathbb{P}_\mu$-almost
surely and $M_\infty(\phi)-M_t(\phi)=o\left(e^{-\frac{\lambda t}{q}}\right)$,
$\mathbb{P}_\mu$-almost surely as $t\to\infty$.
In this case, we have
\begin{equation}\label{id: AA'}
A(q)=\dfrac{q\widetilde A(p)}{\lambda}-\dfrac{q}{\lambda}(M_\infty(\phi)-M_0(\phi)), \quad\mathbb{P}_\mu\mbox{-a.s.}
\end{equation}

$(2)$
 $A_t(q)-A_1(q)$ is in $L^r(\mathbb{P}_\mu)$ and converges  in $L^r(\mathbb{P}_\mu)$ as $t\to\infty$   if and only if  $\widetilde A_t(p)-\widetilde A_1(p)$ is in $L^r(\mathbb{P}_\mu)$ and converges  in $L^r(\mathbb{P}_\mu)$  as $t\to\infty$.  In this case, we have
\begin{equation}\label{id: AA''}
A(q)-A_1(q)=\dfrac{q}{\lambda}(\widetilde A(p)-\widetilde A_1(p))-\dfrac{q}{\lambda}e^{\frac{\lambda}{q}}(M_\infty(\phi)-M_1(\phi))\quad\mbox{ in }L^r(\mathbb{P}_\mu),
\end{equation}
where $A(q)-A_1(q)$ (resp. $\widetilde A(p)-\widetilde A_1(p)$) is the $L^r(\mathbb{P}_\mu)$-limit of $A_t(q)-A_1(q)$ $($resp. $\widetilde A_t(p)-\widetilde A_1(p))$ as $t\to\infty$.
\end{lemma}

\noindent\textbf{Proof:}
The assumption \eqref{LlogL} implies
the uniform integrability of $M_t(\phi), t\in[0,\infty]$. Consequently, $M_t(\phi)$ is bounded on $[0,\infty]$ $\mathbb{P}_\mu$-almost surely.  By the bounded convergence theorem and the stochastic Fubini theorem for martingale measures (c.f.  \cite[Theorem 7.24]{LZ}), for any $T>0$,
\begin{align}\label{identity for A A'}
&A_T(q)=\lim_{l\rightarrow\infty}\int_0^T e^{\frac{\lambda t}{q}} dt\int_t^l e^{-\lambda s}\int_E\phi(x) M(ds ,dx)\\
&= \lim_{l\rightarrow\infty}\int_0^l e^{-\lambda s}\int_E\phi(x) M(ds ,dx)\int_0^{s\wedge T}e^{\frac{\lambda t}{q}} dt\nonumber\\
&= \dfrac{q}{\lambda}\lim_{l\rightarrow\infty}\int_0^l e^{-\lambda s}(e^{\frac{\lambda (s\wedge T)}{q}}-1 )\int_E\phi(x) M(ds ,dx)\nonumber\\
&= \dfrac{q}{\lambda}\int_0^T e^{\frac{-\lambda s}{p}}\int_E\phi(x) M(ds ,dx)+
\dfrac{q}{\lambda}e^{\frac{\lambda T}{q}}\big(M_\infty(\phi)-M_T(\phi)\big)-\dfrac{q}{\lambda}(M_\infty(\phi)-M_0(\phi))\nonumber\\
&= \dfrac{q}{\lambda} \widetilde A_T(p)+\dfrac{q}{\lambda}e^{\frac{\lambda T}{q}}\big(M_\infty(\phi)-M_T(\phi)\big)-\dfrac{q}{\lambda}(M_\infty(\phi)-M_0(\phi)).\nonumber
\end{align}

Note that $A'_t(q):=\dfrac{dA_t(q)}{dt}=e^{\frac{\lambda t}{q}}\big(M_\infty(\phi)-M_t(\phi)\big)$ for almost every $t\in(0,\infty)$.  Therefore \eqref{identity for A A'} can be rewritten as
\begin{equation}\label{equ-At(q)}
-\dfrac{\lambda}{q}A_t(q)+A'_t(q)=(M_\infty(\phi)-M_0(\phi))-\widetilde{A}_t(p), \quad\mbox{a.e. }\, t>0.
\end{equation}
From this we get that
$$
e^{-\frac{\lambda }{q}t}A_t(q)=\frac{q}{\lambda}(M_\infty(\phi)-M_0(\phi))(1-e^{-\frac{\lambda}{q}t})-
\int^t_0e^{-\frac{\lambda}{q}s}\widetilde A_{s}(p)ds, \quad\mbox{a.e. } t>0.
$$
Combining this with \eqref{equ-At(q)}, we  get that for almost every $t\in(0,\infty)$,
\begin{eqnarray}\label{deriv-A}
e^{-\frac{\lambda}{q}t}A'_t(q)
&=&(M_\infty(\phi)-M_0(\phi))(1-e^{-\frac{\lambda}{q}t})-\frac{\lambda}{q}\int^t_0e^{-\frac{\lambda}{q}s}\widetilde A_{s}(p)ds\\
&&
+e^{-\frac{\lambda}{q}t}(M_\infty(\phi)-M_0(\phi))-e^{-\frac{\lambda}{q}t}\widetilde{A}_t(p)\nonumber\\
&=&(M_\infty(\phi)-M_0(\phi))-e^{-\frac{\lambda}{q}t}\widetilde{A}_t(p)-\frac{\lambda}{q}\int^t_0e^{-\frac{\lambda}{q}s}\widetilde A_{s}(p)ds.\nonumber
\end{eqnarray}
Since for a.e. $t>0$, $e^{-\frac{\lambda}{q}t}A'_t(q)=M_\infty(\phi)-M_t(\phi),$ we have for almost all $T,t>0$,
$$
e^{-\frac{\lambda}{q}t}A'_t(q)-e^{-\frac{\lambda}{q}(T+t)}A'_{T+t}(q)=M_{T+t}(\phi)-M_t(\phi).
$$
Using \eqref{deriv-A}, we get that  for almost all $T,t>0$,
\begin{align}\label{derive of A}
&e^{\frac{\lambda}{q}T}\left(M_{T+t}(\phi)-M_T(\phi)\right)\\
&=e^{-\frac{\lambda}{q}t}\widetilde A_{T+t}(p)-\widetilde A_T(p)+\frac{\lambda}{q}\int^t_0e^{-\frac{\lambda}{q}s}\widetilde A_{T+s}ds\nonumber\\
&=(e^{-\frac{\lambda}{q}t}-1)\widetilde A_{T+t}(p)+(\widetilde A_{T+t}(p)-\widetilde A_T(p))+\frac{\lambda}{q}\int^t_0e^{-\frac{\lambda}{q}s}\widetilde A_{T+s}(p)ds .\nonumber
\end{align}
Since
$(1-e^{-\frac{\lambda}{q}t})\widetilde A_{T+t}=\dfrac{\lambda}{q}\int^t_0e^{-\frac{\lambda}{q}s}\widetilde A_{T+t}ds$,
\eqref{derive of A} can be written as: for almost all  $T,t>0$,
\begin{equation}\label{eq: abel}
e^{\frac{\lambda}{q}T}\left(M_{T+t}(\phi)-M_T(\phi)\right)=\widetilde A_{T+t}(p)-\widetilde A_T(p)
-\dfrac{\lambda}{q}\int_0^te^{-\frac{\lambda}{q}s}\left(\widetilde A_{T+t}(p)-\widetilde A_{T+s}(p)\right)ds.
\end{equation}

(1) If $A_T(q)$ converges $\mathbb{P}_\mu$-almost surely as $T\to\infty$ and
$
\lim_{T\to\infty}e^{\frac{\lambda T}{q}}\big(M_\infty(\phi)-M_T(\phi)\big)=0\ \mathbb{P}_\mu\mbox{-a.s.}
$
then by \eqref{identity for A A'},  $\widetilde A_T(p)$ converges  $\mathbb{P}_\mu$-almost surely as $T\to\infty$, and \eqref{id: AA'} follows.  Conversely, if $\widetilde A_T(p)$ converges $\mathbb{P}_\mu$-a.s. as $T\to\infty$, then for any $\varepsilon>0$, there is $\widetilde T(\omega)>0$ such that for $T>\widetilde T(\omega)$ and $t,s\geq 0, \big|\widetilde A_{T+t}(p)-\widetilde A_{T+s}(p)\big|<\varepsilon$.
Using \eqref{eq: abel} and  the right continuity of $M_t(\phi)$, we have for any $t, T>0$,
\begin{eqnarray*}
\Big|e^{\frac{\lambda}{q}T}\left(M_{T+t}(\phi)-M_T(\phi)\right)\Big|&\leq& |\widetilde A_{T+t}(p)-\widetilde A_T(p)|
+\dfrac{\lambda}{q}\int_0^te^{-\frac{\lambda}{q}s}\left|\widetilde A_{T+t}(p)-\widetilde A_{T+s}(p)\right|ds\\
&\leq &\varepsilon+\varepsilon \dfrac{\lambda}{q}\int_0^te^{-\frac{\lambda}{q}s}ds
\leq  2\varepsilon.
\end{eqnarray*}
Letting $t\to\infty$, we get that for $T>\widetilde T$,
$
\Big|e^{\frac{\lambda}{q}T}\left(M_{\infty}(\phi)-M_T(\phi)\right)\Big|\leq 2\varepsilon.
$
Since $\varepsilon>0$ is arbitrary, we have
\[
\mathbb{P}_\mu\left(\lim_{T\to\infty}e^{\frac{\lambda}{q}T}\left(M_{\infty}(\phi)-M_T(\phi)\right)=0\right)=1.
\]
 Thus by \eqref{identity for A A'}, $A_T(q)$ converges $\mathbb{P}_\mu$-almost surely as $T\to\infty$.

(2) Now we  consider the $L^r(\mathbb{P}_\mu)$ convergence.  For any $T\geq 1$, note that
\begin{align*}
&\int_T^{T+1}e^{\frac{\lambda}{q}t}(M_\infty(\phi)-M_t(\phi))dt\nonumber\\
&=\int_T^{T+1}e^{\frac{\lambda}{q}t}(M_\infty(\phi)-M_{T+1}(\phi))dt
+\int_T^{T+1}e^{\frac{\lambda}{q}t}(M_{T+1}(\phi)-M_t(\phi))dt\\
&=\dfrac{q}{\lambda}(M_\infty(\phi)-M_{T+1}(\phi))e^{\frac{\lambda}{q}(T+1)}\left(1-e^{-\lambda/q}\right)+\mathbb P_\mu\left(\int_T^{T+1}e^{\frac{\lambda}{q}t}(M_\infty(\phi)-M_t(\phi))dt\big|\mathcal F_{T+1}\right),\nonumber
\end{align*}
which can be written as
\begin{align}\label{eq: part}
&A_{T+1}(q)-A_T(q)\\
&=\dfrac{q}{\lambda}(M_\infty(\phi)-M_{T+1}(\phi))e^{\frac{\lambda}{q}(T+1)}\left(1-e^{-\lambda/q}\right)+
\mathbb P_\mu\left(A_{T+1}(q)-A_T(q)\big|\mathcal F_{T+1}\right).\nonumber
\end{align}
By Jensen's inequality,
\[
\mathbb P_\mu\left|\mathbb P_\mu\left(A_{T+1}(q)-A_T(q)\big|\mathcal F_{T+1}\right)\right|^r\leq
\mathbb P_\mu\left|A_{T+1}(q)-A_T(q)\right|^r.
\]
If $A_t(q)-A_1(q)$ is in $L^r(\mathbb P_\mu)$ and has
an $L^r(\mathbb P_\mu)$ limit as $t\to\infty$,
then by \eqref{eq: part}, $\{(M_\infty(\phi)-M_{T}(\phi))e^{\frac{\lambda}{q}T}, T\ge 1\}$  is bounded in $L^r(\mathbb P_\mu)$. We obtain from \eqref{identity for A A'} that the  martingale $\{(\widetilde A_t(p)-\widetilde A_1(p)); t\geq 1\}$
is $L^r(\mathbb P_\mu)$ bounded as well.
Thus  the martingale $\widetilde A_t(p)-\widetilde A_1(p)$ has an
$L^r(\mathbb P_\mu)$ limit as $t\to\infty$.

Conversely, if $\widetilde A_t(p)-\widetilde A_1(p)$ is in  $L^r(\mathbb P_\mu)$ and has
an $L^r(\mathbb P_\mu)$ limit $\widetilde A(p)-\widetilde A_1(p)$ as $t\to\infty$,
thanks to \eqref{eq: abel} and Jensen's inequality, for $t\geq 1$,
\begin{align}\label{mart dominated}
&\mathbb P_\mu \left[\big |e^{\frac{\lambda t}{q}}\big(M_\infty(\phi)-M_t(\phi)\big)\big|^r\right]\\
&\lesssim
\mathbb P_\mu|\widetilde{A}(p)-\widetilde{A}_t(p)|^r
+\dfrac{\lambda}{q}\int_0^{\infty} e^{-\frac{\lambda}{q}s}\mathbb P_\mu\left|\widetilde{A}(p)-\widetilde{A}_{t+s}(p)\right|^rds.\nonumber
\end{align}
Applying the dominated convergence theorem to the second term of the right-hand above, we get
\begin{equation}\label{eq:pof di}
\lim_{t\to\infty}\mathbb P_\mu \big |e^{\frac{\lambda t}{q}}\big(M_\infty(\phi)-M_t(\phi)\big)\big|^r=0.
\end{equation}
By Minkowski's inequality, for any $t_1,t_2\geq 1$, in \eqref{identity for A A'}, we deduce that
\begin{align*}
&\displaystyle\|A_{t_1}(q)-A_{t_2}(q)\|_r\\
&\leq \displaystyle\dfrac{q}{\lambda}\left[\|\widetilde{A}_{t_1}(p)-\widetilde{A}_{t_2}(p)\|_r+\|e^{\frac{\lambda t_1}{q}}(M_\infty(\phi)-M_{t_1}(\phi))\|_r+\|e^{\frac{\lambda {t_2}}{q}}(M_\infty(\phi)-M_{t_2}(\phi))\|_r\right].
\end{align*}
Therefore,
$A_t(q)-A_1(q)$ is in $L^r(\mathbb P_\mu)$ and has an $L^r(\mathbb P_\mu)$ limit as $t\to \infty$.\qed
\begin{lemma}\label{le:equi CC'}
Suppose $\mu\in \mathcal{M}^0(E)$ and $\gamma>0$ is a constant.  Define
for $t\ge 0$,
$$
\widetilde{C}_t(\gamma)=\int_0^te^{-\lambda s}s^{\gamma}\int_E\phi(x)M(ds,dx),\quad C_t(\gamma)=\int_0^ts^{\gamma-1}(M_\infty(\phi)-M_s(\phi))ds.
$$
Then $\widetilde{C}_t(\gamma)$ converges $\mathbb{P}_\mu$-almost surely as $t\to\infty$ if and only
if $C_t(\gamma)$ converges and $t^{\gamma-1}(M_\infty(\phi)-M_t(\phi))=o(t^{-1})$, $\mathbb{P}_\mu$-almost surely as $t\to\infty$.
When $\widetilde{C}_t(\gamma)$ (resp. $C_t(\gamma)$) converges as $t\to\infty$, we denote
its limit by $\widetilde C(\gamma)$ (resp. $C(\gamma)$).
 Then we have
\begin{equation}\label{id:CC'}
\gamma\widetilde C(\gamma)=C(\gamma),\quad \mathbb{P}_\mu\mbox{-a.s.}
\end{equation}
\end{lemma}
\noindent\textbf{Proof:}
The proof is similar to that of Lemma \ref{le: equi A A'}.
Similar to \eqref{identity for A A'}, we have for any $T>0$,
\begin{align}\label{identity for CC'}
&C_T(\gamma)=\lim_{l\rightarrow\infty}\int_0^T t^{\gamma-1} dt\int_t^l e^{-\lambda s}\int_E\phi(x) M(ds ,dx)\\
&= \lim_{l\rightarrow\infty}\int_0^l e^{-\lambda s}\int_E\phi(x) M(ds ,dx)\int_0^{s\wedge T}t^{\gamma-1} dt\nonumber\\
&= \frac{1}{\gamma}\lim_{l\rightarrow\infty}\int_0^l e^{-\lambda s}(s\wedge T)^{\gamma}\int_E\phi(x) M(ds ,dx)\nonumber\\
&= \frac{1}{\gamma}\int_0^T e^{-\lambda s}s^{\gamma}\int_E\phi(x) M(ds ,dx)+
\frac{1}{\gamma} T^{\gamma}\big(M_\infty(\phi)-M_T(\phi)\big)\nonumber\\
&= \frac{1}{\gamma}\widetilde C_T(\gamma)+\frac{1}{\gamma} T^{\gamma}\big[M_\infty(\phi)-M_T(\phi)\big].\nonumber
\end{align}
If
$C_t(\gamma)$ converges and $t^{\gamma-1}(M_\infty(\phi)-M_t(\phi))=o(t^{-1})$ as $t\to\infty$ $\mathbb{P}_\mu$-almost surely.  Using \eqref{identity for CC'}, we get
$\widetilde C_T(\gamma)$ converges $\mathbb{P}_\mu$-almost surely and \eqref{id:CC'} holds.
We now deduce the almost sure convergence of $C_t(\gamma)$ and $t^{\gamma}\big[M_\infty(\phi)-M_t(\phi)\big]$ from the a.s. convergence of $\widetilde C_t(\gamma)$. From \eqref{identity for CC'}, we get
\begin{eqnarray*}\label{id:c1}
T^{\gamma}\big[M_\infty(\phi)-M_T(\phi)\big]=\gamma C_T(\gamma)-\widetilde C_T(\gamma)
\end{eqnarray*}
and
\begin{eqnarray*}\label{id:c2}
T^{\gamma}\big[M_\infty(\phi)-M_{T+t}(\phi)\big]=\dfrac{T^{\gamma}}{(T+t)^{\gamma}}\left[\gamma C_{T+t}(\gamma)-\widetilde C_{T+t}(\gamma)\right].
\end{eqnarray*}
It follows from the two displays above that
\begin{eqnarray}\label{id: C3}
\qquad T^{\gamma}\big[M_{T+t}(\phi)-M_{T}(\phi)\big]=\gamma T^{\gamma}\left[\dfrac{C_{T}(\gamma)}{T^\gamma}-\dfrac{C_{T+t}(\gamma)}{(T+t)^{\gamma}}\right]+\left[\dfrac{T^{\gamma}}{(T+t)^{\gamma}}\widetilde C_{T+t}(\gamma)-\widetilde C_{T}(\gamma)\right].
\end{eqnarray}
Noticing that $t^\gamma(M_\infty(\phi)-M_t(\phi))=t \dfrac{dC_t(\gamma)}{dt}=tC'_t(\gamma)$
for a.e. $t>0$, \eqref{identity for CC'} can be written as
\begin{equation*}\label{id: CC' ode}
\gamma C_t(\gamma)-t C'_t(\gamma)=\widetilde{C}_t(\gamma), \quad\mbox{ a.e. }t>0.
\end{equation*}
From this we get that for any $t,T>0$,
\begin{eqnarray}\label{sol}
\dfrac{C_{T}(\gamma)}{T^\gamma}-\dfrac{C_{T+t}(\gamma)}{(T+t)^{\gamma}}=\int_T^{T+t}\dfrac{\widetilde C_s(\gamma)}{s^{\gamma+1}}ds.
\end{eqnarray}
Simple calculations yield
\begin{eqnarray}\label{id:2}
\dfrac{T^{\gamma}}{(T+t)^{\gamma}}\widetilde C_{T+t}(\gamma)-\widetilde C_{T}(\gamma)= \widetilde C_{T+t}(\gamma)-\widetilde C_{T}(\gamma)-\gamma T^{\gamma}\int_T^{T+t}\dfrac{\widetilde C_{T+t}(\gamma)}{s^{\gamma+1}}ds.
\end{eqnarray}
Combining \eqref{id: C3}, \eqref{sol} and \eqref{id:2}, we get
\begin{equation}\label{id: abel2}
T^\gamma\left[M_{T+t}(\phi)-M_{T}(\phi)\right]=\widetilde{C}_{T+t}(\gamma)-\widetilde{C}_{T}(\gamma)-\gamma T^\gamma\int_T^{T+t}\dfrac{\widetilde{C}_{T+t}(\gamma)-\widetilde{C}_{s}(\gamma)}{s^{\gamma+1}}ds.
\end{equation}
If $\widetilde{C}_t(\gamma)$ converges  $\mathbb{P}_\mu$-almost surely as $t\to\infty$,
by using the argument in the proof of Lemma \ref{le: equi A A'}(1) after \eqref{eq: abel}, we get
\begin{equation}
\lim_{t\to\infty}t^\gamma\left[M_{\infty}(\phi)-M_t(\phi)\right]=0,\quad\mathbb P_\mu\mbox{-a.s.}
\end{equation}
Therefore,
it follows from \eqref{identity for CC'} that $C_{t}(\gamma)$ converges $\mathbb{P}_{\mu}$-almost surely
and \eqref{id:CC'} holds as well. \qed
\begin{remark}\label{remark for general CC'}
Suppose that $L(ds,dx)$ is a random measure on $[0,\infty)\times E$
such that, as $t\to\infty$,
$$
   L_t:=\int_0^t e^{-\lambda s}\int_E\phi(x)L(ds,dx)\to L_\infty,\quad \mathbb{P}_\mu\mbox{-a.s.}
$$
where $L_\infty$ is a finite random variable.
Using arguments similar to those in the proof of Lemma \ref{le:equi CC'}, we can show
that for any $\gamma>0$, $\int_0^T e^{-\lambda s}s^\gamma\int_E\phi(x)L(ds,dx)$
converges $\mathbb{P}_\mu$-almost surely as $T\to \infty$ if and only if
$$
 \int_0^T t^{\gamma-1}dt\int_t^\infty e^{-\lambda s}\int_E\phi(x)L(ds,dx)=\int^T_0t^{\gamma-1}(L_\infty-L_t)dt
$$
converges and $L_\infty-L_T=o(T^{-\gamma})$,
$\mathbb{P}_\mu$-almost surely as $T\to \infty$.
\end{remark}

\begin{lemma}\label{lem:Lpconv}
  Assume that Assumptions \ref{asp:H1}-\ref{asp:H2}  and \eqref{LlogL} hold. Let $1\leq a<p\leq 2$.
\begin{itemize}
\item[(1)] If \eqref{finite of p moment} holds, then for any $\mu\in\mathcal M^0(E)$, $(\widetilde A_t(a)-\widetilde A_1(a))$ is in $L^p(\mathbb P_\mu)$  and converges in $L^p(\mathbb P_\mu)$ and therefore
$\mathbb{P}_\mu$-almost surely as $t\to\infty$.
\item[(2)] Suppose that for some $\mu\in\mathcal M^0(E)$, $(\widetilde A_t(a)-\widetilde A_1(a))$
  is in $L^p(\mathbb P_\mu)$  and converges in $L^p(\mathbb P_\mu)$
 as $t\to \infty$,
then it must converge
$\mathbb{P}_\mu$-almost surely  as $t\to \infty$
and \eqref{finite of p moment} holds.
\end{itemize}
\end{lemma}
\textbf{Proof:}
(1) Suppose condition \eqref{finite of p moment} holds.  From the definition \eqref{def: ano mart} of  $\widetilde A_t(a)$ and \eqref{M=J+C}, we only need to consider the convergences of
$$
\int_1^{t} e^{-\frac{\lambda}{a} s}\int_E \phi(x)S^J(ds, dx)\qquad \mbox{and}\qquad
\int_1^{t} e^{-\frac{\lambda}{a} s}\int_E \phi(x)S^C(ds, dx)
$$
as $t\to\infty$.
Recall that definitions of $S^{(1,\infty)}_t$  and $S^{(2,\infty)}_t$ given in \eqref{martingale-int2} and \eqref{martingale-int2(2)} with $\rho=\infty$.
For the ``small jump" part, by the Burkholder-Davis-Gundy inequality, we have
\begin{eqnarray*}
\mathbb P_\mu\left[(\sup_{t\geq 0}S^{(1,\infty)}_t(e^{-\frac{\lambda}{a} \cdot}\phi))^2\right]
&\lesssim &\mathbb P_\mu\left(\int_0^\infty\int_{\mathcal{M}(E_\partial)}F^2_{e^{-\frac{\lambda}{a}\cdot}\phi}(s,\upsilon)\widehat{N}^{(1,\infty)}(ds,d\upsilon)\right)\\
&=&\mathbb P_{\mu}\int_{0}^{\infty} ds\int_E{X}_s(dx)
\int_0^{1}F^2_{e^{-\frac{\lambda}{a}\cdot}\phi}(s, r\phi(x)^{-1}\delta_x)\pi^\phi(x,dr)\\
&=&\int_{0}^{\infty}e^{-\frac{2\lambda}{a} s}\langle P_s^\beta\left(\int_0^{1} r^2\pi^\phi(\cdot,dr)\right),\mu\rangle ds.
\end{eqnarray*}
Thanks to \eqref{assum: moment finite}, we have
\begin{align}\label{ine-small}
 &\mathbb P_\mu\left[\left(\sup_{t\geq 0}S^{(1,\infty)}_t(e^{-\frac{\lambda}{a} \cdot}\phi)\right)^2\right]
\lesssim \int_{0}^{\infty}e^{-\frac{2\lambda}{a} s} \langle P_s^\beta\phi,\mu\rangle ds\\
&=\dfrac{a}{(2-a)\lambda}\int_E\phi(y)\mu(dy)<\infty,\nonumber
\end{align}
where in the last equality, we used the fact that $e^{-\lambda s}P_s^\beta\phi=\phi$ for $s\geq 0$.
Applying the Burkholder-Davis-Gundy inequality to $S^{(2,\infty)}_t(e^{-\frac{\lambda}{a} \cdot}\phi))$, we obtain that for $1<p\leq 2$,
\begin{align*}
&\mathbb P_\mu\left[\Big|\sup_{t\geq 1}\big(S^{(2,\infty)}_t(e^{-\frac{\lambda}{a} \cdot}\phi)-S^{(2,\infty)}_1(e^{-\frac{\lambda}{a} \cdot}\phi)\big)\Big|^p\right]\lesssim \mathbb P_\mu\bigg[\sum_{s\in [1,\infty)\cap
J^{(2,\infty)}}F_{e^{-\frac{\lambda}{a}\cdot}\phi}(s, \Delta \overline X_s)^2\bigg]^{\frac{p}{2}}\\
&\leq  \mathbb P_\mu\bigg[\sum_{s\in [1,\infty)\cap J^{(2,\infty)}}F_{e^{-\frac{\lambda}{a}\cdot}\phi}(s, \Delta \overline X_s)^p\bigg]\\
&=\mathbb P_{\mu}\int_1^\infty ds\int_E {X}_s(dx)
\int_1^{\infty}F^p_{e^{-\frac{\lambda}{a}\cdot}\phi}(s, r\phi(x)^{-1}\delta_x)\pi^\phi(x,dr)\\
&=\int_1^{\infty}e^{-\frac{p\lambda}{a} s}\langle P_s^\beta\left(\int_{1}^\infty r^p\pi^\phi(\cdot,dr)\right),\mu\rangle ds.
\end{align*}
Set $h(x):=\int_{1}^\infty r^p\pi^\phi(x,dr)$. Condition \eqref{finite of p moment} says that $h\in L_1^+(\nu)$. Note that $p>a$. If $\mu\in\mathcal M^0(E)$, by Assumption \ref{asp:H2},
\[
\int_1^{\infty}e^{-\frac{p\lambda}{a} s}\langle P_s^\beta\left(\int_{1}^\infty r^p\pi^\phi(\cdot,dr)\right),\mu\rangle ds\lesssim
 \mu(\phi)\nu(h)\Big(1+\sup_{t>1,x\in E}|C_{t,x,h}|\Big)\int_1^{\infty}e^{-\left(\frac{p}{a}-1\right)\lambda s}ds<\infty.
\]
Therefore,
\begin{eqnarray}\label{ine-large}
\mathbb P_\mu\left[\Big|\sup_{t\geq 1}\big(S^{(2,\infty)}_t(e^{-\frac{\lambda}{a}\cdot}\phi)-S^{(2,\infty)}_1(e^{-\frac{\lambda}{a} \cdot}\phi)\big)\Big|^p\right]<\infty.
\end{eqnarray}
Combining \eqref{ine-small} and \eqref{ine-large}, we get
\begin{eqnarray}\label{ine-jump}
\mathbb P_\mu\left[\Big|\sup_{t\geq 1}\int_1^{t} e^{-\frac{\lambda}{a} s}\int_E \phi(x)S^J(ds, dx)\Big|^p\right]<\infty.
\end{eqnarray}
We also have
\begin{align}\label{ine-con}
&\sup_{t\geq 0}\mathbb P_\mu\left[\Big(\int_0^{t} e^{-\frac{\lambda}{a} s}\int_E \phi(x)S^C(d s, dx)\Big)^2\right]
=\mathbb P_\mu\int_0^\infty e^{-\frac{2}{a}\lambda s}\langle \alpha\phi^2, {X}_s\rangle ds\\
&=\int_0^\infty e^{-\frac{2}{a}\lambda s}ds\int_EP_s^{\beta}(\alpha\phi^2)(y)\mu(dy)
\leq \dfrac{a\|\alpha\phi\|_\infty}{(2-a)\lambda}\langle\phi,\mu\rangle <\infty.\nonumber
\end{align}
Consequently, by \eqref{ine-jump} and \eqref{ine-con}, $\sup_{t\geq 1}\mathbb P_\mu\left(|\widetilde A_t(a)-\widetilde A_1(a)|^p\right)<\infty$.
Thus $\widetilde A_t(a)-\widetilde A_1(a)$  is in $L^p(\mathbb P_\mu)$ and
converges in $L^p(\mathbb P_\mu)$ and $\mathbb P_\mu$-almost
surely as $t\to\infty$.

(2) Suppose that, for some $\mu\in \mathcal M^0(E)$,
$\widetilde A_t(a)-\widetilde A_1(a)\to \widetilde A(a)-\widetilde A_1(a)$ in $L^p(\mathbb P_\mu)$ as $t\to \infty$. Then $\mathbb P_\mu\left(|\widetilde A(a)-\widetilde A_1(a)|^p\right)<\infty$.
 Since $\widetilde A_t(a)$ is a $\mathbb P_\mu$-martingale,
it must converge $\mathbb P_\mu$-almost surely as $t\to\infty$.
By Jensen's inequality, for any $t>1$,
\[
\mathbb P_\mu( |\widetilde A_t(a)-\widetilde A_1(a)|^p)=\mathbb P_\mu \left(\left|\mathbb P_\mu\left(\widetilde A(a)-\widetilde A_1(a)\big|\mathcal F_t\right)\right|^p\right)\leq \mathbb P_\mu( |\widetilde A(a)-\widetilde A_1(a)|^p)<\infty.
\]
We have shown in \eqref{ine-small} and \eqref{ine-con} that
\[
\mathbb P_\mu\left[\left(S^{(1,\infty)}_t(e^{-\frac{\lambda}{a} \cdot}\phi)\right)^2\right]<\infty
\qquad \mbox{and}\qquad
\mathbb P_\mu\left[\Big(\int_0^{t} e^{-\frac{\lambda}{a} s}\int_E \phi(x)S^C(d s, dx)\Big)^2\right]<\infty.
\]
Therefore, by the definition of $\widetilde A_t(a)$ given in \eqref{def: ano mart}, we have that for any $t\geq 0$,
\begin{equation}\label{big jump finite}
\mathbb P_\mu\left(\left|S^{(2,\infty)}_t(e^{-\frac{\lambda}{a} \cdot}\phi)-S^{(2,\infty)}_1(e^{-\frac{\lambda}{a} \cdot}\phi)\right|^p\right)<\infty.
\end{equation}
Note that it
follows from \eqref{mart dominated} that $\mathbb P_\mu(\left|M_\infty(\phi)-M_1(\phi)\right|^p)<\infty$ when
$\mathbb P_\mu (|\widetilde A(a)-\widetilde A_1(a)|^p)<\infty$.  Therefore,
\[
\mathbb P_\mu\left[M_\infty(\phi)^p\Big|\mathcal F_1\right]\lesssim \mathbb P_\mu\left[|M_\infty(\phi)-M_1(\phi)|^p\Big|\mathcal F_1\right]+M_1(\phi)^p<\infty.
\]
Since $\{M_t(\phi);t\geq 1\}$ is a martingale with respect to $(\mathcal F_t)_{t\geq 1}$ under $\mathbb P_\mu(\cdot|\mathcal F_1)$,
we have almost surely
\[
\mathbb P_\mu\left[\sup_{t\geq 1}M_t(\phi)^p\Big|\mathcal F_1\right]\lesssim \mathbb P_\mu\left[M_\infty(\phi)^p\Big|\mathcal F_1\right]<\infty.
\]
Thus for the compensator $\widehat{N}^{(2,\infty)}$ of the ``big jumps", we have $\mathbb P_\mu$-almost surely
\begin{align}\label{comp finit}
&\mathbb P_\mu\left[\Big(\int_1^t\int_{\mathcal{M}(E_\partial)}F_{e^{-\frac{\lambda}{a}\cdot}\phi}(s,\upsilon)\widehat{N}^{(2,\infty)}(ds,d\upsilon)\Big)^p\Big|\mathcal F_1\right]\\
&=\mathbb P_\mu\left[\Big(\int_1^te^{-\frac{\lambda}{a} s}ds\int_{E}{X}_s(dx)
\int_1^\infty r \pi^\phi(x,dr)\Big)^p\Big|\mathcal F_1\right]\nonumber\\
&\lesssim  \mathbb P_\mu\left[\Big(\int_1^te^{\frac{\lambda}{a^*} s}M_s(\phi)ds\Big)^p\Big|\mathcal F_1\right]
\leq f^p_a(t)\mathbb P_\mu\left(\sup_{s\leq t}M_s(\phi)^p\Big|\mathcal F_1\right)\nonumber\\
&\lesssim  f^p_a(t)\mathbb P_\mu\left(M_{t}(\phi)^p\Big|\mathcal F_1\right)
\leq  f^p_a(t)\mathbb P_\mu\left(M_{\infty}(\phi)^p\Big|\mathcal F_1\right)<\infty,\quad t\geq 1,\nonumber
\end{align}
where $\dfrac{1}{a^*}+\dfrac{1}{a}=1$ and
$
f_a(t)=\begin{cases}\dfrac{a^*}{\lambda}e^{\frac{\lambda}{a^*}t},& a>1\\
t,& a=1,
\end{cases}
$ and in the first inequality we used \eqref{assum: moment finite}.
It follows  from \eqref{big jump finite} and \eqref{comp finit} that
for any $t>1$,
\[
\mathbb P_\mu\left[\Big(\int_1^t\int_{\mathcal{M}(E_\partial)}F_{e^{-\frac{\lambda}{a}\cdot}\phi}(s,\upsilon)N^{(2,\infty)}(ds,d\upsilon)\Big)^p\Big|\mathcal F_1\right]<\infty.
\]
Since $p>1$, $\{X_t\}$ is Markov and $F_{e^{-\lambda\cdot}\phi}(s,\upsilon)\geq 0$, for any $t>0$,
\begin{align*}
 \infty&>\mathbb P_{X_1}\left[\Big(\int_0^{t}\int_{\mathcal{M}(E_\partial)}F_{e^{-\frac{\lambda}{a}\cdot}\phi}(s+1,\upsilon)N^{(2,\infty)}(ds, d\upsilon)\Big)^p\right]\\
&= \mathbb P_{X_1}\left[\Big(\sum_{s\in (0,t]\bigcap J^{(2,+\infty)}}F_{e^{-\frac{\lambda}{a}\cdot}\phi}(s+1,\Delta \overline X_s)\Big)^p\right]\\
&\geq  \mathbb P_{X_1}\left[\sum_{s\in (0,t]\bigcap J^{(2,+\infty)}}\Big(F_{e^{-\frac{\lambda}{a}\cdot}\phi}(s+1,\Delta \overline X_s)\Big)^p\right]\\
&= \mathbb P_{X_1}\Big(\int_0^t e^{-\frac{p\lambda}{a} (s+1)}ds\int_E {X}_s(dx)
\int_1^\infty r^p \pi^\phi(x,dr)\Big)\\
&= \int^t_0 e^{-\frac{p\lambda}{a} (s+1)}ds\int_EP_s^\beta\Big(\int_1^\infty r^p \pi^\phi(\cdot,dr)\Big)(y)X_1(dy),
\end{align*}
which implies that for any $T\geq 1$,
\[
\int^{1+T}_T e^{-\frac{p\lambda}{a} s}ds\int_EP_s^\beta h(y)X_1(dy)<\infty,
\]
where $h(x)=\int_1^\infty r^p \pi^\phi(x,dr)$ as before.  Since $\nu$ is a probability measure, $h\wedge L\in L_1^+(\nu)$ for any $L>0$.
Let $c_t=\sup_{x\in E, f\in L_1^+(\nu) }\left|C_{t,x,f}\right|$ be as in Assumption \ref{asp:H2}. Then $\lim_{t\to\infty}c_t=0$. Choose
$T>0$ such that when $t>T$, $c_t\leq 1/2$. Applying Assumption \ref{asp:H2} again, we get  that for $t>T$,
\[
e^{\lambda t}\phi(y)\nu(h\wedge L)\leq 2P_t^\beta h(y).
\]
Integrating both sides of the above inequality with respect to
$e^{-\frac{p}{a}\lambda t}dtX_1(dy)$
on $[T,T+1]\times E$, we obtain
\begin{align*}
&\dfrac{\nu(h\wedge L)X_1(\phi)}{(\frac{p}{a}-1)\lambda}e^{-(\frac{p}{a}-1)\lambda T}\left(1-e^{-(\frac{p}{a}-1)\lambda}\right)\leq \int_T^{T+1}e^{-\frac{p}{a}\lambda t}dt\int_E\nu(h\wedge L)e^{\lambda t}\phi(y)X_1(dy)\\
&\leq 2\int^{1+T}_T e^{-\frac{p}{a}\lambda t}dt\int_EP_s^\beta h(y)X_1(dy)<\infty.
\end{align*}
Since the last term in the above inequality does not depend on $L$, letting
 $L\to\infty$ we get
\[
\nu(h)=\int_E \nu(dx)\int_1^\infty r^p \pi^\phi(x,dr)<\infty.
\]
The proof is now complete.\qed

\begin{lemma}\label{le: upper}
Suppose \eqref{assm: unif upp} holds and $T_0$ is the constant in \eqref{assm: unif upp}. Then there is a constant $K>0$  such that,  for any $\mu\in\mathcal M^0(E)$, and any $t,n, a, b>0$ satisfying $0<b<a$ and $e^{bn}>T_0$,
\[
\widetilde{\mathbb P}_\mu(\langle\phi, X_t\rangle>e^{a n})\leq 3\langle\phi,\mu\rangle e^{\lambda t-an}+Ke^{\lambda t-(a-b)n}+Kt\int_E\nu(dy)\int_{e^{bn}}^\infty r\pi^\phi(y,dr).
\]
\end{lemma}
\textbf{Proof:}
It follows from the spine decomposition that
\begin{align}\label{domi-sums}
&\widetilde{\mathbb P}_\mu(\langle\phi, X_t\rangle>e^{a n})
=\mathbb Q_\mu\Big(\langle\phi,X_t\rangle+\langle\phi,Z^{C}_t\rangle+\langle\phi,Z^{J}_t\rangle>e^{an}
\Big)\nonumber\\
&\leq  \mathbb Q_\mu\left(\langle\phi,X_t\rangle>\dfrac{1}{3}e^{an}\right)+\mathbb Q_\mu\left(\langle\phi,Z^{C}_t\rangle>\dfrac{1}{3}e^{an}\right)
+\mathbb Q_\mu\left(\langle\phi,Z^{J}_t\rangle>\dfrac{1}{3}e^{an}\right)\\
&\leq \dfrac{3\mathbb Q_\mu\langle\phi,X_t\rangle}{e^{a n}}+\dfrac{3\mathbb Q_\mu\langle\phi,Z^{C}_t\rangle}{e^{a n}}+\mathbb Q_\mu\left(\langle\phi,Z^{J}_t\rangle>1/3e^{an}\right).\nonumber
\end{align}
Noting that $\mathbb Q_\mu\langle\phi,X_t\rangle=\mathbb P_\mu\langle\phi,X_t\rangle=e^{\lambda t}\langle\phi,\mu\rangle$, we get
$$
\dfrac{3\mathbb Q_\mu\langle\phi,X_t\rangle}{e^{a n}}\le 3\langle\phi,\mu\rangle e^{\lambda t-an}.
$$
Since
$
\mathbb Q_\mu\langle\phi,Z^{C}_t\rangle=\widetilde\Pi_{\phi\mu}\left[\int_0^t\alpha(\xi_s)\phi(\xi_s)e^{\lambda(t-s)}ds\right]\leq \dfrac{2\|\alpha\phi\|_\infty e^{\lambda t}}{\lambda},
$
it follows that
\begin{eqnarray}\label{iq: 2ed}
\dfrac{3\mathbb Q_\mu\langle\phi,Z^{C}_t\rangle}{e^{a n}}\leq \dfrac{6\|\alpha\phi\|_\infty e^{\lambda t}e^{-an}}{\lambda}.
\end{eqnarray}
Let $m_\sigma= X^{J,\sigma}_0(E)$ for $\sigma\in \mathcal D^{J}$.
From the construction of $Z_t^{J}$, we can estimate  the third term in \eqref{domi-sums} as follows:
\begin{eqnarray*}
\mathbb Q_\mu\left(\langle\phi,Z^{J}_t\rangle>\dfrac{1}{3}e^{an}\right)\leq \mathbb Q_\mu\left(\sum_{\stackrel{\sigma \in \mathcal D^J \cap [0, t]}{
m_\sigma\phi(\xi_\sigma)>e^{bn}}}\langle\phi, X^{J,\sigma}_{t-\sigma}\rangle>\dfrac{1}{6}e^{an}\right)+\mathbb Q_\mu\left(\sum_{\stackrel{\sigma \in \mathcal D^J \cap [0, t]}{m_\sigma\phi(\xi_\sigma)\leq e^{bn}}}\langle\phi, X^{J,\sigma}_{t-\sigma}\rangle>\dfrac{1}{6}e^{an}\right).
\end{eqnarray*}
By the Markov inequality,
\begin{align}\label{iq:2ed2}
& \mathbb Q_\mu\left(\sum_{\stackrel{\sigma \in \mathcal D^J \cap [0, t]}{
m_\sigma\phi(\xi_\sigma)\leq e^{bn}}}\langle\phi, X^{J,\sigma}_{t-\sigma}\rangle>\dfrac{1}{6}e^{an}\right)\leq 6e^{-an}\mathbb Q_\mu\left(\sum_{\stackrel{\sigma \in \mathcal D^J \cap [0, t]}{
m_\sigma\phi(\xi_\sigma)\leq e^{bn}}}\langle\phi, X^{J,\sigma}_{t-\sigma}\rangle\right)\\
&= 6e^{-an}\widetilde\Pi_{\phi\mu}\left[\int_0^te^{\lambda(t-s)}\phi^{-1}(\xi_s)ds\int_0^{e^{bn}}r^2\pi^{\phi}(\xi_s,dr)\right]\nonumber\\
&=6e^{-an}\widetilde\Pi_{\phi\mu}\int_0^te^{\lambda(t-s)}\phi^{-1}(\xi_s)ds\left(\int_0^{1}r^2\pi^\phi(\xi_s,dr)+\int_1^{e^{bn}}
r^2\pi^\phi(\xi_s,dr)\right)\nonumber\\
&\leq 6e^{-(a-b)n}\widetilde\Pi_{\phi\mu}\int_0^te^{\lambda(t-s)}\phi^{-1}(\xi_s)ds\int_0^{\infty}(r\wedge r^2)\pi^\phi(\xi_s,dr).\nonumber
\end{align}
Thus by  \eqref{assum: moment finite} and \eqref{iq:2ed2},
\begin{equation}\label{iq:2ed2'}
\mathbb Q_\mu\left(\sum_{\stackrel{\sigma \in \mathcal D^J\cap [0, t]}{
m_\sigma\phi(\xi_\sigma)\leq e^{bn}}}\langle\phi, X^{J,\sigma}_{t-\sigma}\rangle>\dfrac{1}{6}e^{an}\right)
 \leq \dfrac{6C}{\lambda}e^{\lambda t-(a-b)n}:=Ae^{\lambda t-(a-b)n},
\end{equation}
where $A=\dfrac{6C}{\lambda}$
and $C$ is the constant in \eqref{assum: moment finite}.  It is obvious that $A$ is independent of $\mu, t,a$ and  $b$.
When the event
$\left\{\sum_{\stackrel{\sigma \in \mathcal D^J \cap [0, t]}{m_\sigma\phi(\xi_{\sigma})>e^{bn}}}
\langle\phi, X^{J,\sigma}_{t-\sigma}\rangle>\dfrac{1}{6}e^{an}
\right\}$ occurs, $\sharp\{\sigma \in \mathcal D^J \cap [0, t]; m_\sigma\phi(\xi_{\sigma})>e^{bn}\}\ge 1$.
Thus
\begin{align*}
&\mathbb Q_\mu\left(\sum_{\stackrel{\sigma \in \mathcal D^J \cap [0, t]}{m_\sigma\phi(\xi_{\sigma})>e^{bn}}}\langle\phi, X^{J,\sigma}_{t-\sigma}\rangle>\dfrac{1}{6}e^{an}\right)\leq \mathbb Q_\mu\left(\sum_{\stackrel{\sigma \in \mathcal D^J \cap [0, t]}{
m_\sigma\phi(\xi_{\sigma})>e^{bn}}}1\geq 1\right)\\
&\leq \mathbb Q_\mu\left(\sum_{\stackrel{\sigma \in \mathcal D^J \cap [0, t]}{m_\sigma\phi(\xi_{\sigma})>e^{bn}}}1\right)
=\widetilde{\Pi}_{\phi\mu}\left[\int_0^t ds\int_{e^{bn}}^\infty r\pi^{\phi}(\xi_s,dr)\right].
\end{align*}
When $e^{bn}\geq T_0$, we have
\[
\widetilde{\Pi}_{\phi\mu}\left[\int_0^t ds\int_{e^{bn}}^\infty r\pi^\phi(\xi_s,dr)\right]
\leq B\widetilde{\Pi}_{\phi\mu}\left[\int_0^t \phi(\xi_s)ds \int_E\nu(dy)\int_{e^{bn}}^\infty r\pi^\phi(y,dr)\right].
\]
Therefore,
\begin{eqnarray}\label{iq: 3rd}
\mathbb Q_\mu\left(\sum_{\stackrel{\sigma \in \mathcal D^J \cap [0, t]}{m_\sigma\phi(\xi_{\sigma})>e^{bn}}}\langle\phi, X^{\mathrm m,\sigma}_{t-\sigma}\rangle>\dfrac{1}{6}e^{an}\right)\leq Bt\|\phi\|_\infty \int_E\nu(dy)\int_{e^{bn}}^\infty r \pi^\phi(y,dr).
\end{eqnarray}
Put $K=\dfrac{6\|\alpha\phi\|_\infty}{\lambda}+A+B\|\phi\|_\infty$, which is independent of  $\mu, t,a$ and $b$. Combining \eqref{domi-sums}, \eqref{iq: 2ed}, \eqref{iq:2ed2'} and \eqref{iq: 3rd}, we obtain
\[
\widetilde{\mathbb P}_\mu(\langle\phi, X_t\rangle>e^{a n})\leq 3\langle\phi,\mu\rangle e^{\lambda t-an}+Ke^{\lambda t-(a-b)n}+Kt\int_E\nu(dy)\int_{e^{bn}}^\infty r\pi^\phi(y,dr).
\]
\qed

\subsection{Proofs of Main Results}
In this subsection, we give the proofs of our main results.

\noindent\textbf{Proof of Theorem \ref{theorem:Lpconv}: }
(1) Suppose \eqref{finite of p moment} holds. Using  Lemma \ref{lem:Lpconv}(1)
with $1<a< p\le 2$,
$(\widetilde A_t(a)-\widetilde A_1(a))$ converges in $L^p(\mathbb P_\mu)$ and $\mathbb P_\mu$-almost surely as $t\to\infty$.
Then by Lemma \ref{le: equi A A'}(2),
$(A_t(a^*)-A_1(a^*))$ converges in $L^p(\mathbb P_\mu)$ as $t\to\infty$.

(2) Suppose that for some $\mu\in\mathcal M^0(E)$, $(A_t(a^*)-A_1(a^*))$  converges in $L^p(\mathbb P_\mu)$ as $t\to \infty$.
By Lemma \ref{le: equi A A'}(2), $(\widetilde A_t(a)-\widetilde A_1(a))$ converges
in $L^p(\mathbb P_\mu)$ as $t\to\infty$.
Applying Lemma \ref{lem:Lpconv}(2), we get that it converges
$\mathbb P_\mu$-almost surely as $t\to\infty$
and \eqref{finite of p moment} holds.
It now follows from Lemma \ref{le: equi A A'}(1) that $(A_t(a^*)-A_1(a^*))$  converges  $\mathbb P_\mu$-almost surely as $t\to \infty$.

(3) According to (1),  $(\widetilde A_t(a^*)-\widetilde A_1(a^*))$ converges in
$L^p(\mathbb P_\mu)$ as $t\to\infty$.
Repeating the argument leading to \eqref{eq:pof di}, we get
$$
\lim_{t\to\infty}\mathbb P_\mu \big |e^{\frac{\lambda t}{a^*}}\big(M_\infty(\phi)-M_t(\phi)\big)\big|^p=0.
$$
Thus the assertion of (3) holds.

(4) This  is the result  of Lemma \ref{lem:Lpconv}(2) with $a=1$. \qed

\noindent\textbf{Proof of Theorem \ref{thm: p moment convergence rate}: }
(1) Suppose  \eqref{finite of p moment} holds  and $\mu\in \mathcal{M}^0(E)$.
 Since $(\widetilde A_t(p))$ is a martingale, it converges
 $\mathbb{P}_\mu$-almost surely as $t\to \infty$ if it is uniformly integrable.
Note that
$$
\widetilde A_t(p)=\int_0^1 e^{\frac{-\lambda s}{p}}\int_E\phi(x)M(ds,dx)+\int_1^t e^{\frac{-\lambda s}{p}}\int_E\phi(x)M(ds,dx),\qquad t\geq 1.
$$
We only need to consider the convergence of
$$
\int_1^t e^{\frac{-\lambda s}{p}}\int_E\phi(x)M(ds,dx)
$$
as $t\to\infty$.

For  the ``small jumps" part, we have, for $t>0$,
\begin{align*}
& \mathbb P_\mu\left[\Big(\int_1^te^{-\frac{\lambda}{p} s}\int_E \phi(x) S^{(1,p)}(ds, dx)\Big)^2\right]\\
&\lesssim \int_{1}^{t}e^{\frac{-2\lambda}{p} s} ds \int_E P_s^\beta\Big(\int_0^{e^{\frac{\lambda}{p}s}} r^2\pi^\phi(\cdot,dr)\Big)(y)\mu(dy)\\
&\lesssim \mu(\phi)\int_{1}^{\infty}e^{(1-\frac{2}{p})\lambda s} ds\int_E \nu(dx)\int_0^{e^{\frac{\lambda}{p}s}} r^2\pi^\phi(x,dr)\\
&\leq \int_{1}^{\infty}e^{\lambda s(1-\frac{2}{p})} ds\int_E \nu(dx)\int_0^{1} r^2\pi^\phi(x,dr)+\int_{1}^{\infty}e^{\lambda s(1-\frac{2}{p})} ds\int_E \nu(dx)\int_1^{e^{\frac{\lambda}{p}s}} r^2\pi^\phi(x,dr)\\
&=I+II,
\end{align*}
where in the second inequality we used Assumption \ref{asp:H2}.
Since $p<2$, we have $I<\infty$.
When \eqref{finite of p moment} holds, by Fubini's theorem, we get
\begin{eqnarray*}
II&\leq &\int_E \nu(dx)\int_1^{\infty} r^2\pi^\phi(x,dr)\int_{\frac{p}{\lambda}\ln r}^{\infty}e^{\lambda s(1-\frac{2}{p})} ds\\
&\lesssim&\int_E \nu(dx)\int_1^{\infty} r^p\pi^\phi(x,dr)<\infty.
\end{eqnarray*}
It follows that
\begin{equation}\label{small-jump}
\sup_{t>1}{\mathbb P}_\mu\left[\Big(\int_1^te^{-\frac{\lambda}{p} s}\int_E \phi(x) S^{(1,p)}(ds, dx)\Big)^2\right]<\infty.
\end{equation}

For  the ``big jumps" part, using Assumption \ref{asp:H2} and Fubini's theorem again, we get
\begin{align}\label{big-jump}
& \mathbb P_\mu\sup_{t>1}\Big|\int_1^{t}e^{-\frac{\lambda}{p} s}\int_E \phi(x) S^{(2,p)}(ds, dx)\Big|\\
&\leq  2\mathbb P_\mu\int_1^\infty ds\int_E {X}_s(dx)
\int_{e^{\lambda s/p}}^{\infty}F_{e^{-\frac{\lambda}{p}\cdot}\phi}(s, r\phi(x)^{-1}\delta_x)\pi^\phi(x,dr)\nonumber\\
&\lesssim \int_1^\infty e^{\frac{\lambda}{q} s}ds\int_E\nu(dx)\int_{e^{\frac{\lambda}{p}s}}^\infty r\pi^\phi(x,dr)\nonumber\\
&\lesssim   \int_E\nu(dx)\int_{1}^\infty r^p \pi^\phi(x,dr)<\infty.\nonumber
\end{align}
For the continuum part, we have the following estimates:
\begin{eqnarray*}
\sup_{t>1}\mathbb P_\mu\left[\left(\int_1^{t}e^{-\frac{\lambda}{p} s}\int_E \phi(x)S^C(ds, dx)\right)^2\right]
&=&\mathbb P_\mu\int_1^\infty e^{-\frac{2\lambda}{p} s}ds\int_E\alpha(x)\phi(x)^2 {X}_s(dx)\\
&\lesssim & \int_1^\infty e^{-\lambda s(2/p-1)}ds\int_E\alpha(x)\phi(x)^2\nu(dx).
\end{eqnarray*}
Since $p<2$,
\begin{equation}\label{con-part}
\sup_{t>1}\mathbb P_\mu\left[\left(\int_1^{t}e^{-\frac{\lambda}{p} s}\int_E \phi(x)S^C(ds, dx)\right)^2\right]<\infty.
\end{equation}
Combining \eqref{small-jump}, \eqref{big-jump} and \eqref{con-part}, we obtain that $\int_1^t e^{\frac{-\lambda s}{p}}\int_E\phi(x)M(ds,dx)$ is uniformly integrable.
 Thus $\widetilde{A}_t(p)$ converges $\mathbb{P}_\mu$-almost surely as $t\to\infty$.
 By Lemma \ref{le: equi A A'}, $A_t(q)$ converges $\mathbb{P}_\mu$-almost surely as $t\to\infty$ and
\[
M_\infty(\phi)-M_t(\phi)=o\left(e^{-\frac{\lambda}{q}t}\right),\qquad \mathbb P_\mu\mbox{-a.s. }\, {\rm  as }\,\, t\to\infty.
\]

(2) Now we suppose
 \begin{equation}\label{inf p mom}
\int_E\nu(dy)\int_1^\infty r^p \pi^\phi(y, dr)=\infty
\end{equation}
and  $\mu\in \mathcal{M}^0(E)$.
By Assumption \ref{asp:H2}, there is $t_0>0$ such that for any $f\in L_1^+(\nu)$ and $t>t_0$,
\begin{equation}\label{IU>}
P_t^\beta f(x) \ge \dfrac{1}{2}e^{\lambda t}\phi(x)\nu(f),\quad x\in E.
\end{equation}
Without loss of generality, we assume $t_0=1/2$.
Set $\rho_t=e^{\lambda t/q}(M_\infty(\phi)-M_t(\phi))$, $t>0$.
For any $n\in\mathbb N$ and $1/2\leq t\leq 1$, note that
\begin{align*}
&\Delta\rho_{n+t}=-e^{-\lambda (n+t)/p}\Delta \overline X_{n+t}(\phi),
\end{align*}
and thus  $\Delta \overline X_{n+t}(\phi)>2e^{\lambda n/p}$ implies that
$$
\left|\rho_{n+t}\right|>e^{-\lambda/p}\quad \mbox{or}\quad \left|\rho_{(n+t)-}\right|>e^{-\lambda/p}.
$$
Define
\[
B_n=\left\{\sum_{\frac{1}{2}\leq t<1}1_{\{\triangle \overline X_{n+t}(\phi)>2e^{\lambda n/p}\}}>0\right\}.
\]
 Then we have $\{B_n,{\rm i.o.}\}$ implies
 $\rho_t=o(1)$ does not hold as $t\to\infty$ .
 Therefore, we only need to prove $\mathbb P_\mu(B_n,{\rm i.o.})>0$.
 If we can prove that
 \begin{equation}\label{BC-toprove}
\sum_{n=1}^\infty\mathbb P_\mu\left(B_n\Big|\mathcal F_n\right)=\infty,\quad \mbox{ a.s.\, on }\, \{M_\infty(\phi)>0\},
\end{equation}
then by the  second conditional Borel-Cantelli  lemma (see, \cite[Theorem 5.3.2]{Durrett}),
$$
\mathbb P_\mu(B_n,{\rm i.o.})=\mathbb{P}_\mu\Big(\sum_{n=1}^\infty\mathbb P_\mu\left(B_n\Big|\mathcal F_n\right)=\infty\Big)\geq \mathbb P_\mu\Big(M_\infty(\phi)>0\Big)>0.
$$
Therefore, we only need to prove \eqref{BC-toprove}.

To prove \eqref{BC-toprove}, we will estimate the probability ${\bf P}(Y>0)$ for the non-negative random variable $Y:=
\sum_{n+\frac{1}{2}\leq t<n+1}1_{\{\triangle \overline X_t(\phi)>2e^{\lambda n/p}\}}$ defined on some probability space with probability ${\bf P}$.
Our basic idea is to use  the inequality  ${\bf P}(Y>0)\geq \frac{({\bf P}Y)^2}{{\bf P}Y^2} $.
However, $Y$ may not have second moment.  Thus we consider $\sum_{n+1/2\leq t<n+1}1_{\{\triangle \overline X_t(\phi)>2e^{\lambda n/p}\}\bigcap C^A_n(t)}$ for
 some appropriate events $C^A_n(t)$.  We will prove \eqref{BC-toprove} in 4 steps.

{\bf Step 1.}\quad We first prove that
\begin{equation}\label{inf-exp}
\sum_{n=1}^\infty\mathbb P_\mu\Big(\sum_{n+\frac{1}{2}\leq t<n+1}1_{\{\triangle\overline X_t(\phi)>2e^{\lambda n/p}\}}\Big|\mathcal F_n\Big)=\infty \quad {\rm a.s.\,\, on}\quad \{M_\infty(\phi)>0\}.
\end{equation}
Using \eqref{IU>} with $t_0=1/2$,  we have
\begin{align}\label{ine-2}
&\mathbb P_\mu\Big(\sum_{n+\frac{1}{2}\leq t<n+1}1_{\{\triangle \overline X_t(\phi)>2e^{\lambda n/p}\}}\Big|\mathcal F_n\Big)=\Big\langle\int_{\frac{1}{2}}^1ds P_s^\beta\Big(\int_{2e^{\lambda n/p}}^\infty \pi^\phi(\cdot, dr)\Big), X_n\Big\rangle\\
&\gtrsim\langle\phi, X_n\rangle \int_E\nu(dx)\int_{2e^{\lambda n/p}}^\infty \pi^\phi(x, dr).\nonumber
\end{align}
  Therefore,
\begin{equation}\label{ine-l}
\sum_{n=1}^\infty\mathbb P_\mu\Big(\sum_{n+\frac{1}{2}\leq t<n+1}1_{\{\triangle \overline X_t(\phi)>2e^{\lambda n/p}\}}\Big|\mathcal F_n\Big)
\gtrsim\sum_{n=1}^\infty\langle\phi, X_n\rangle\int_E\nu(dx)\int_{2e^{\lambda n/p}}^\infty \pi^\phi(x, dr).
\end{equation}
Since $\lim_{n\rightarrow\infty}e^{-\lambda n}\langle\phi, X_n\rangle =M_\infty(\phi)$ almost surely, on the event $\{M_\infty(\phi)>0\}$, the convergence of the right-hand side of \eqref{ine-l} is equivalent to the convergence of the series
$$
\sum_{n=1}^\infty e^{\lambda n} \int_E\nu(dx)\int_{2e^{\lambda n/p}}^\infty \pi^\phi(x, dr).
$$
Since $e^{\lambda s}$ is increasing and $\int_{2e^{\lambda s/p}}^\infty \pi^\phi(x, dr)$ is decreasing in $s$, we have
\begin{align*}
&\sum_{n=1}^\infty e^{\lambda n} \int_E\nu(dx)\int_{2e^{\lambda n/p}}^\infty \pi^\phi(x, dr)\\
&\geq \int_0^\infty e^{\lambda s} ds\int_E\nu(dx)\int_{2e^{\lambda/p}e^{\lambda s/p}}^\infty \pi^\phi(x, dr)\\
&= \int_E\nu(dx)\int_{2e^{\lambda/p}}^\infty \pi^\phi(x, dr)\int_0^{\frac{p}{\lambda}\ln r-\frac{p}{\lambda}\ln (2e^{\lambda/p})} e^{\lambda s} ds\\
&\gtrsim \int_E\nu(dx)\int_{2e^{\lambda/p}}^\infty r^p\pi^\phi(x, dr).
\end{align*}
Therefore \eqref{inf-exp} follows from \eqref{inf p mom}.

 {\bf Step 2.}\quad Suppose $A>0$ is an arbitrary fixed constant.  For any integer $n\geq 1$, define
\begin{equation}\label{def:C}
C^{A}_n(t)=\left\{\dfrac{\langle\int_{2e^{\lambda n/p}}\pi^\phi(\cdot, dr), \overline X_{t-}\rangle}{e^{\lambda n}\langle\int_{2e^{\lambda n/p}}\pi^\phi(\cdot, dr),\nu\rangle}<L\right\}, \quad n+\frac{1}{2}\leq t<n+1,
\end{equation}
and
\begin{equation}
\label{def:C>} C^{>}_n(t)=\left\{\dfrac{\langle \int_{2e^{\lambda n/p}}\pi^\phi(\cdot, dr), \overline X_{t-}\rangle}{e^{\lambda n}\langle\int_{2e^{\lambda n/p}}\pi^\phi(\cdot, dr),\nu\rangle}>L/2\right\},\quad \frac{1}{2}\leq t<1,
\end{equation}
where $L$ is chosen large enough so that for
any $\mu\in\mathcal M^0(E)$ satisfying $\langle\phi,\mu\rangle<A$,
\begin{equation}\label{def: L}
\mathbb P_\mu\left((C^A_n(t))^c\right)<\dfrac{1}{4},\quad n\geq 1,\, n+\frac{1}{2}\le t<n,
\end{equation}
and for any $\mu\in\mathcal M^0(E)$ satisfying $\langle\phi,\mu\rangle<Ae^{\lambda n}$,
\begin{equation}\label{def: L/2}
\mathbb P_\mu\left(C^{>}_n(t)\right)<\dfrac{1}{2},\quad n\geq 1,\,t\in[\frac{1}{2},1].
\end{equation}
The existence of such an $L$ is guaranteed by
\begin{align*}
&\mathbb P_\mu\left(\dfrac{\langle\int_{2e^{\lambda n/p}}\pi^\phi(\cdot, dr), \overline X_{t-}\rangle}{e^{\lambda n}\langle\int_{2e^{\lambda n/p}}\pi^\phi(\cdot, dr),\nu\rangle}>L\right)\leq
\dfrac{\mathbb P_\mu\left(\langle\int_{2e^{\lambda n/p}}\pi^\phi(\cdot, dr), \overline X_{t-}\rangle\right)}{L e^{\lambda n}\langle\int_{2e^{\lambda n/p}}\pi^\phi(\cdot, dr),\nu\rangle}\\
&\leq \dfrac{(1+c_t)e^{\lambda t}\langle\phi, \mu\rangle}{Le^{\lambda n}}
\leq\dfrac{A e^\lambda(1+c_t)}{L},
\end{align*}
if $n+\dfrac{1}{2}\leq t<n+1$ and $\langle\phi,\mu\rangle<A$, or if $\dfrac{1}{2}\leq t<1$ and $\langle\phi,\mu\rangle<Ae^{\lambda n}$.
The first inequality above is the Markov inequality, and $c_t$ is the
quantity  in Assumption \ref{asp:H2} which is bounded for $t>1/2$.
Thus $L$ can be chosen large enough to assure both \eqref{def: L} and \eqref{def: L/2} hold.
In this step, we will prove that there is $N\in\mathbb N$ such that when $n>N$,
$\mathbb P_\mu$-almost surely on $\{M_n(\phi)\leq A\}$,
\begin{align}\label{withA-withoutA}
&\mathbb P_\mu\Big(\sum_{n+1/2\leq t<n+1}1_{\{\triangle \overline X_t(\phi)>2e^{\lambda n/p}\}\bigcap C^A_n(t)}\Big|\mathcal F_n\Big)\\
&\geq \dfrac{1}{4}\mathbb P_\mu\left(\sum_{n+1/2\leq t<n+1}1_{\{\triangle \overline X_t(\phi)>2e^{\lambda n/p}\}}\Big|\mathcal F_n\right).\nonumber
\end{align}
We divide $X_n$ into $[e^{\lambda n}]$ disjoint parts each with value $[e^{\lambda n}]^{-1}X_n$.
For $i=1,2,\ldots ,[e^{\lambda n}]$, let $(X_s^{(i)},\, 0\leq s<1)$ be the superprocess with the $i$-th part as its initial mass.
By the branching property of  superprocesses, $X_s^{(i)},\, i=1,2,\ldots ,[e^{\lambda  n}]$, are independent and identically distributed as $\mathbb P_{[e^{\lambda n}]^{-1}X_n}$
under $\mathbb P_\mu\left(\cdot\Big|\mathcal F_n\right)=\mathbb P_{X_n}(\cdot)$.
Thus for any $i=1,2,\ldots ,[e^{\lambda  n}]$,
\begin{align*}
  &(C_n^A(s))^c=\left\{\dfrac{\langle\int_{2e^{\lambda n/p}}\pi^\phi(\cdot, dr), \overline X_{s-}\rangle}{e^{\lambda n}\langle\int_{2e^{\lambda n/p}}\pi^\phi(\cdot, dr), \nu \rangle}>L\right\}\\
  &\subset \left\{\dfrac{\langle\int_{2e^{\lambda n/p}}\pi^\phi(\cdot, dr), \overline X^{(i)}_{s-}\rangle}{e^{\lambda n}\langle\int_{2e^{\lambda n/p}}\pi^\phi(\cdot, dr), \nu\rangle}>L/2\right\}
  \bigcup\left\{\dfrac{\sum_{j\neq i}\langle\int_{2e^{\lambda n/p}}\pi^\phi(\cdot, dr), \overline X^{(j)}_{s-}\rangle}
  {e^{\lambda n}\langle\int_{2e^{\lambda n/p}}\pi^\phi(\cdot, dr),\nu\rangle}>L/2\right\}\\
  &:=C_n^{(i)}(s)\bigcup C_n^{(\neq i)}(s),
  \quad n+\frac{1}{2}\leq s<n+1, n\geq 1.
\end{align*}
Consider the conditional expectation:
\begin{align*}
&E_n:=\mathbb P_\mu\Big(\sum_{n+\frac{1}{2}\leq t<n+1}1_{\{\triangle \overline X_t(\phi)>2e^{\lambda n/p}\}\setminus C^A_n(t)}\Big|\mathcal F_n\Big)\\
&=\mathbb P_{X_n}\left(\int_{\frac{1}{2}}^1 1_{ (C^A_n(s))^c}ds\int_E\overline X_{s-}(dx)\int_{2e^{\frac{\lambda}{p}n}}^\infty \pi^{\phi}(x, dr)\right)\\
&=\sum_{i=1}^{[e^{\lambda n}]}\mathbb P_{X_n}\left(\int_{\frac{1}{2}}^1 1_{ (C^A_n(s))^c}ds\int_E\overline X^{(i)}_{s-}(dx)\int_{2e^{\frac{\lambda}{p}n}}^\infty \pi^{\phi}(x, dr)\right)\\
&\leq \sum_{i=1}^{[e^{\lambda n}]}\mathbb P_{X_n}\left(\int_{\frac{1}{2}}^1 (1_{C_n^{(i)}(s)}+1_{C_n^{(\neq i)}(s)})ds\int_E\overline X^{(i)}_{s-}(dx)\int_{2e^{\frac{\lambda}{p}n}}^\infty \pi^{\phi}(x, dr)\right)\\
&\leq \sum_{i=1}^{[e^{\lambda n}]}\mathbb P_{X_n}\left(\int_{\frac{1}{2}}^1 1_{C_n^{(i)}(s)} ds\int_E\overline X^{(i)}_{s-}(dx)\int_{2e^{\frac{\lambda}{p}n}}^\infty \pi^{\phi}(x, dr)\right)\\
&\qquad\qquad+\sum_{i=1}^{[e^{\lambda n}]}\mathbb P_{X_n}\left(\int_{\frac{1}{2}}^1 1_{C_n^{(\neq i)}(s)}ds\int_E\overline X^{(i)}_{s-}(dx)\int_{e^{\frac{\lambda}{p}n}}^\infty \pi^{\phi}(x, dr)\right)\\
&:= I_n^{(1)}+I_n^{(2)}.
\end{align*}
Since $X^{(i)}$ and $X^{(\neq i)}=\sum_{j\neq i}X^{(j)}$ are independent,
\begin{align*}
1_{\{M_n(\phi)\leq A\}}I_n^{(2)}&=\sum_{i=1}^{[e^{\lambda n}]}I_{\{M_n(\phi)\leq A\}}\mathbb P_{X_n}\left(\int_{\frac{1}{2}}^1 1_{C_n^{(\neq i)}(s)}ds\int_E\overline X^{(i)}_{s-}(dx)\int_{e^{\frac{\lambda}{p}n}}^\infty \pi^{\phi}(x, dr)\right)\\
&=\sum_{i=1}^{[e^{\lambda n}]}\int_{\frac{1}{2}}^1 ds 1_{\{M_n(\phi)\leq A\}}\mathbb P_{X_n}\big(C_n^{(\neq i)}(s)\big)\mathbb P_{X_n}\left(\int_E\overline X^{(i)}_{s-}(dx)\int_{e^{\frac{\lambda}{p}n}}^\infty \pi^{\phi}(x, dr)\right)\\
&\leq\sum_{i=1}^{[e^{\lambda n}]}\int_{\frac{1}{2}}^1 ds 1_{\{M_n(\phi)\leq A\}}\mathbb P_{X_n}\big(C_n^{>}(s)\big)\mathbb P_{X_n}\left(\int_E\overline X^{(i)}_{s-}(dx)\int_{e^{\frac{\lambda}{p}n}}^\infty \pi^{\phi}(x,dr)\right).
\end{align*}
On the event $\{M_n(\phi)\leq A\}$, we have $\langle\phi, X_n\rangle\le A e^{\lambda n}$.
Therefore,
\begin{equation}\label{In(2)}
1_{\{M_n(\phi)\leq A\}}I_n^{(2)}\leq 1_{\{M_n(\phi)\leq A\}}\dfrac{1}{2}\mathbb P_\mu\left(\sum_{n+\frac{1}{2}\leq t<n+1}1_{\{\triangle \overline X_t(\phi)>2e^{\lambda n/p}\}}\big|\mathcal F_n\right).
\end{equation}
As for $I_n^{(1)}$, when $2e^{\frac{\lambda}{p}n}>T_0$, by assumption \eqref{assm: unif upp},
\[
\int_E\overline X_{s-}^{(i)}(dx)\int_{2e^{\frac{\lambda}{p}n}}^\infty \pi^\phi(x,dr)\leq B \langle\phi,\overline X_{s-}^{(i)}\rangle \int_E\nu(dx)\int_{2e^{\frac{\lambda}{p}n}}^\infty \pi^\phi(x,dr),\quad s\in [\frac{1}{2}, 1).
\]
Therefore,
\[
C_n^{(i)}(s)\subset\{\langle\phi,\overline X_{s-}^{(i)}\rangle\geq \frac{L}{2B}e^{\lambda n}\},\quad s\in [\frac{1}{2}, 1)
\]
and
\begin{align*}
& I_n^{(1)}= \sum_{i=1}^{[e^{\lambda n}]}\mathbb P_{X_n}\left(\int_{\frac{1}{2}}^1 1_{C_n^{(i)}(s)} ds
\int_E\overline X^{(i)}_{s-}(dx)\int_{2e^{\frac{\lambda}{p}n}}^\infty \pi^{\phi}(x,dr)\right)\\
&\leq  Be^{\lambda n}\int_E\nu(dx)\int_{2e^{\frac{\lambda}{p}n}}^\infty \pi^\phi(x,dr)
\mathbb P_{[e^{\lambda n}]^{-1}X_n}\int_{\frac{1}{2}}^1 \langle\phi, \overline X_{s-}\rangle 1_{\{\langle\phi, \overline X_{s-}\rangle
\geq \frac{L}{2B}e^{\lambda n}\}} ds\nonumber \\
&=  Be^{\lambda n}\int_E\nu(dx)\int_{2e^{\frac{\lambda}{p}n}}^\infty \pi^\phi(x,dr)
\int_{\frac{1}{2}}^1 \mathbb P_{[e^{\lambda n}]^{-1}X_n}\langle\phi, X_s\rangle 1_{\{\langle\phi,X_s\rangle\geq \frac{L}{2B}e^{\lambda n}\}} ds\nonumber \\
&=  Be^{\lambda n}[e^{\lambda n}]^{-1}X_n(\phi)\int_E\nu(dx)\int_{2e^{\frac{\lambda}{p}n}}^\infty \pi^\phi(x,dr)\int_{\frac{1}{2}}^1 e^{\lambda s} ds \widetilde{\mathbb P}_{[e^{\lambda n}]^{-1}X_n}\left(\langle\phi,X_s\rangle\geq \frac{L}{2B}e^{\lambda n}\right).\nonumber
\end{align*}
Note that we may choose $L$ large enough that $\frac{L}{2B}\geq 1$.
From Lemma \ref{le: upper}, for any $0<b<\lambda$, and any $1/2<s<1$, on the set $\{M_n(\phi)\leq A\}$, there is a constant $\widetilde K>0$ such that almost surely
\begin{align*}
&\widetilde{\mathbb P}_{[e^{\lambda n}]^{-1}X_n}\left(\langle\phi, X_s\rangle\geq \frac{L}{2B}e^{\lambda n}\right)\leq \widetilde{\mathbb P}_{[e^{\lambda n}]^{-1}X_n}\left(\langle\phi, X_s\rangle\geq e^{\lambda n}\right)\\
&\leq
3[e^{\lambda n}]^{-1}\langle\phi,X_n\rangle e^{\lambda s-\lambda n}+Ke^{\lambda s-(\lambda-b)n}+Ks\int_E\nu(dy)\int_{e^{bn}}^\infty r\pi^\phi(y,dr)\\
&\leq  \widetilde Ke^{-(\lambda-b)n}+\widetilde K\int_E\nu(dy)\int_{e^{bn}}^\infty r\pi^\phi(y,dr).
\end{align*}
Therefore,
\begin{align*}
&1_{\{M_n(\phi)\leq A\}}I_n^{(1)}\\
&\lesssim 1_{\{M_n(\phi)\leq A\}}X_n(\phi)\left(\int_E\nu(dx)\int_{2e^{\frac{\lambda}{p}n}}^\infty \pi^\phi(x,dr)\right)\left[e^{-(\lambda-b)n}+\int_E\nu(dy)\int_{e^{bn}}^\infty r\pi^\phi(y,dr)\right]\\
&\lesssim  1_{\{M_n(\phi)\leq A\}}\mathbb P_\mu\Big(\sum_{n+\frac{1}{2}\leq t<n+1}1_{\{\triangle \overline X_t(\phi)>2e^{\lambda n/p}\}}\Big|\mathcal F_n\Big)\left[e^{-(\lambda-b)n}+\int_E\nu(dy)\int_{e^{bn}}^\infty r\pi^\phi(y,dr)\right],
\end{align*}
where the last inequality follows from \eqref{ine-2}.  Since
$\lim_{n\to\infty}e^{-(\lambda-b)n}+\int_E\nu(dy)\int_{e^{bn}}^\infty rn^\phi(y,dr)=0$, we can choose $N>0$ such that when $n\geq N$, we have
 $e^{bn}>T_0$ and
 \begin{align}\label{In(1)}
 1_{\{M_n(\phi)\leq A\}}I_n^{(1)}\leq \dfrac{1}{4}1_{\{M_n(\phi)\leq A\}}\mathbb P_\mu\Big(\sum_{n+\frac{1}{2}\leq t<n+1}1_{\{\triangle \overline X_t(\phi)>2e^{\lambda n/p}\}}\Big|\mathcal F_n\Big),\quad \mathbb{P}_\mu\mbox{-a.s.}
 \end{align}
 Combining \eqref{In(2)} and \eqref{In(1)}, we get, when $n>N$, on $\{M_n(\phi)\leq A\}$,
 $$
 E_n\leq I_n^{(1)}+I_n^{(2)}\leq \dfrac{3}{4}1_{\{M_n(\phi)\leq A\}}\mathbb P_\mu\Big(\sum_{n+\frac{1}{2}\leq t<n+1}1_{\{\triangle \overline X_t(\phi)>2e^{\lambda n/p}\}}\Big|\mathcal F_n\Big),\quad \mathbb{P}_\mu\mbox{-a.s.}
 $$
Therefore, when $n>N$, on $\{M_n(\phi)\leq A\}$,
\begin{align*}
& \mathbb P_\mu\Big(\sum_{n+1/2\leq t<n+1}1_{\{\triangle \overline X_t(\phi)>2e^{\lambda n/p}\}\bigcap C^A_n(t)}\Big|\mathcal F_n\Big)\\
&\geq \mathbb P_\mu\Big(\sum_{n+1/2\leq t<n+1}1_{\{\triangle \overline X_t(\phi)>2e^{\lambda n/p}\}}\Big|\mathcal F_n\Big)-E_n\\
&\geq  \dfrac{1}{4}\mathbb P_\mu\left(\sum_{n+1/2\leq t<n+1}1_{\{\triangle \overline X_t(\phi)>2e^{\lambda n/p}\}}\big|\mathcal F_n\right),\quad \mathbb{P}_\mu\mbox{-a.s.}
\end{align*}
This proves \eqref{withA-withoutA}.

 {\bf Step 3.}\quad In this step, we prove that, on $\{M_\infty(\phi)>0, \sup_{n}M_n(\phi)\leq A\}$,
\begin{equation}\label{step2-sum}
\sum_{n=1}^\infty\mathbb P_\mu\left(\sum_{n+1/2\leq t<n+1}1_{\{\triangle \overline X_t(\phi)>2e^{\frac{\lambda n}{p}}\}\bigcap C^A_n(t)}>0\Big|\mathcal F_n\right)=\infty,\quad \mathbb{P}_\mu\mbox{-a.s.}
\end{equation}

Let $N$ be a number large enough  so that \eqref{withA-withoutA} almost surely holds on $\{M_n(\phi)\leq A\}$ for any $n\ge N$. Then on the event $\{M_\infty(\phi)>0,\sup_{n}M_n(\phi)\leq A\}$,
\begin{align*}
&\sum_{n=N}^m 1_{\{M_n(\phi)\leq A\}}\mathbb P_\mu\Big(\sum_{n+1/2\leq t<n+1}1_{\{\triangle \overline X_t(\phi)>2e^{\lambda n/p}\}\bigcap C^A_n(t)}\Big|\mathcal F_n\Big)\\
&\geq  \dfrac{1}{4}\sum_{n=N}^m1_{\{M_n(\phi)\leq A\}}\mathbb P_\mu\left(\sum_{n+1/2\leq t<n+1}1_{\{\triangle \overline X_t(\phi)>2e^{\lambda n/p}\}}\Big|\mathcal F_n\right)\\
&=  \dfrac{1}{4}\sum_{n=N}^m\mathbb P_\mu\left(\sum_{n+1/2\leq t<n+1}1_{\{\triangle \overline X_t(\phi)>2e^{\lambda n/p}\}}\Big|\mathcal F_n\right),\quad\mathbb{P}_\mu\mbox{-a.s.}
\end{align*}
By \eqref{inf-exp}, letting $m\to\infty$ in the display above, we get that,
on $\{\sup_{n>1}M_n(\phi)\leq A, M_\infty (\phi)>0\}$,
$$
\sum_{n=1}^\infty\mathbb P_\mu\Big(\sum_{n+1/2\leq t<n+1}1_{\{\triangle \overline X_t(\phi)>2e^{\lambda n/p}\}\bigcap C^A_n(t)}\Big|\mathcal F_n\Big)=\infty,\quad \mathbb P_\mu\mbox{-a.s.}
$$

 Let $\widetilde C_n(t):=\Big\{\triangle \overline X_t(\phi)>2e^{\frac{\lambda n}{p}}\Big\}\bigcap C^A_n(t)$.  Now we consider the second moments:
\begin{align*}
&\mathbb P_\mu\left[\Big(\sum_{n+1/2\leq t<n+1}1_{\{\triangle \overline X_t(\phi)>2e^{\frac{\lambda n}{p}}\}\bigcap C^A_n(t)}\Big)^2\Big|\mathcal F_n\right]\\
&=2\mathbb P_\mu\Big[\sum_{n+1/2\leq t_1<t_2<n+1}1_{\widetilde C_n(t_1)}1_{\widetilde C_n(t_2)}\Big|\mathcal F_n\Big]+\mathbb P_\mu\Big[\sum_{n+1/2\leq t<n+1}1_{\big\{\triangle \overline X_t(\phi)>2e^{\frac{\lambda n}{p}}\big\}\bigcap C^A_n(t)}\Big|\mathcal F_n\Big].
\end{align*}
Define $I_n^{(3)}:=2\mathbb P_\mu\Big[\sum_{n+1/2\leq t_1<t_2<n+1}1_{\widetilde C_n(t_1)}1_{\widetilde C_n(t_2)}\Big|\mathcal F_n\Big]$, then
\begin{align*}
&I_n^{(3)}=2\mathbb P_\mu\Big[\sum_{n+1/2\leq t_1<n+1}1_{\widetilde C_n(t_1)}\mathbb P_\mu \left(\sum_{t_1<t_2<n+1}1_{\widetilde C_n(t_2)}\Big|\mathcal F_{t_1}\right)\Big|\mathcal F_n\Big]\\
&\leq  2\mathbb P_\mu\Big[\sum_{n+1/2\leq t_1<n+1}1_{\widetilde C_n(t_1)}\mathbb P_{X_{t_1}}\left(\int_{t_1}^{n+1}1_{C^A_n(s)}ds \int_E \overline X_{s-}(dx)\int_{2e^{\lambda n/p}}^\infty \pi^{\phi}(x,dr)\right)\Big|\mathcal F_n\Big]\\
&\leq  2Le^{\lambda n}\int_E\nu(dx)\int_{2e^{\lambda n/p}}^\infty \pi^\phi(x, dr)\cdot\mathbb P_\mu\Big[\sum_{n+1/2\leq t_1<n+1}1_{\widetilde C_n(t_1)}\Big|\mathcal F_n\Big]\\
&\lesssim e^{\lambda n}\langle\phi,X_n\rangle^{-1} \mathbb P^2_\mu\Big[\sum_{n+1/2\leq t_1<n+1}1_{\big\{\triangle \overline X_{t_1}(\phi)>2e^{\frac{\lambda n}{p}}\big\}}\Big|\mathcal F_n\Big],
\end{align*}
where the last inequality comes from \eqref{ine-2}. Consequently
\begin{align}\label{est:second}
&\mathbb P_\mu\left[\Big(\sum_{n+1/2\leq t<n+1}1_{\{\triangle \overline X_t(\phi)>2e^{\frac{\lambda n}{p}}\}\bigcap C^A_n(t)}\Big)^2\Big|\mathcal F_n\right]\\
&\lesssim
\frac1{M_n(\phi)}
\mathbb P^2_\mu\Big[\sum_{n+1/2\leq t<n+1}1_{\big\{\triangle \overline X_{t}(\phi)>2e^{\frac{\lambda n}{p}}\big\}}\Big|\mathcal F_n\Big]+\mathbb P_\mu\Big[\sum_{n+1/2\leq t<n+1}1_{\big\{\triangle \overline X_t(\phi)>2e^{\frac{\lambda n}{p}}\big\}\bigcap C^A_n(t)}\Big|\mathcal F_n\Big].\nonumber
\end{align}
Now by the Cauchy-Schwarz inequality, when $n>N$, on the event $\{M_\infty(\phi)>0,\sup_{n}M_n(\phi)\leq A\}$,
\begin{align*}
&\mathbb P_\mu\left(\sum_{n+1/2\leq t<n+1}1_{\{\triangle \overline X_t(\phi)>2e^{\frac{\lambda n}{p}}\}\bigcap C^A_n(t)}>0\Big|\mathcal F_n\right)\\
&\geq \dfrac{\mathbb P^2_\mu\left(\sum_{n+1/2\leq t<n+1}1_{\{\triangle \overline X_t(\phi)>2e^{\frac{\lambda n}{p}}\}\bigcap C^A_n(t)}\Big|\mathcal F_n\right)}{\mathbb P_\mu\left(\Big(\sum_{n+1/2\leq t<n+1}1_{\{\triangle \overline X_t(\phi)>2e^{\frac{\lambda n}{p}}\}\bigcap C^A_n(t)}\Big)^2\Big|\mathcal F_n\right)}\\
&\gtrsim  \dfrac{\dfrac{1}{16}\mathbb P^2_\mu\left(\sum_{n+1/2\leq t<n+1}1_{\{\triangle \overline X_t(\phi)>2e^{\frac{\lambda n}{p}}\}}\Big|\mathcal F_n\right)}{M_n(\phi)^{-1} \mathbb P^2_\mu\Big[\sum_{n+1/2\leq t_1<n+1}1_{\big\{\triangle \overline X_{t_1}(\phi)>2e^{\frac{\lambda n}{p}}\big\}}\Big|\mathcal F_n\Big]+\mathbb P_\mu\Big[\sum_{n+1/2\leq t_1<n+1}1_{\big\{\triangle \overline X_{t_1}(\phi)>2e^{\frac{\lambda n}{p}}\big\}}\Big|\mathcal F_n\Big]}\\
&\gtrsim  \dfrac{\mathbb P_\mu\left(\sum_{n+1/2\leq t<n+1}1_{\{\triangle \overline X_t(\phi)>2e^{\frac{\lambda n}{p}}\}}\Big|\mathcal F_n\right)}{M_n(\phi)^{-1} \mathbb P_\mu\Big[\sum_{n+1/2\leq t_1<n+1}1_{\big\{\triangle \overline X_{t_1}(\phi)>2e^{\frac{\lambda n}{p}}\big\}}\Big|\mathcal F_n\Big]+1}\\
&\gtrsim  M_n(\phi)\bigwedge \mathbb P_\mu\left(\sum_{n+1/2\leq t<n+1}1_{\{\triangle \overline X_t(\phi)>2e^{\frac{\lambda n}{p}}\}}\Big|\mathcal F_n\right),\quad \mathbb P_\mu\mbox{-a.s.}
\end{align*}
where the second inequality comes from \eqref{withA-withoutA} and \eqref{est:second}, and the last inequality comes from the fact that
$\frac{x}{y+1}\geq \frac{x}{2(y\wedge 1)}\gtrsim \frac{x}{y}\wedge x$ for any $x,y>0$.  Since we
are working
on $\{M_\infty(\phi)>0\}$ and
we have proved  \eqref{inf-exp}, \eqref{step2-sum} follows from the inequalities above.

{\bf Step 4.}\quad By the  second conditional Borel-Cantelli lemma (see, \cite[Theorem 5.3.2]{Durrett}), \eqref{step2-sum} implies that
\[
\mathbb P_\mu\left(B_n \, {\rm i.o.} \, \Big|M_\infty(\phi)>0,\sup_{n}M_n(\phi)\leq A\right)=1.
\]
Note that $\sup_{n\geq 1}M_n(\phi)<\infty$ $\mathbb P_\mu$-almost surely.  The above equation  holds for any constant $A>0$.  Letting $A\to\infty$, we get \eqref{BC-toprove}. Consequently,
$$
\limsup_{t\to\infty}\left|e^{\lambda t/q}\left(M_\infty(\phi)-M_t(\phi)\right)\right|\geq e^{-\lambda/p}
$$
with positive probability. The proof is complete.
\qed

\bigskip

For $\gamma>0$, define
\begin{equation}\label{def-f}
f(s)=e^{\lambda s} s^{-\gamma}, s>0.
\end{equation}
Direct computation shows that
\[
f'(s)=f(s)(\lambda-\gamma s^{-1}).
\]
Thus when $s>\gamma/\lambda$, $f(s)$ is a strictly increasing function. If $g$ is the inverse function of $f$ on $(\gamma/\lambda,\infty)$, then
\[
(g(r))'=\dfrac{1}{r(\lambda-\gamma g(r)^{-1})}.
\]
It is obvious that
\[
\lim_{r\to\infty}g(r)=\infty.
\]
Therefore, there is a constant $R>\gamma/\lambda$ such that for $r>R$,
\begin{equation}\label{eq: R}
\dfrac{1}{\lambda r}\leq (g(r))'\leq \dfrac{2}{\lambda r}.
\end{equation}
Consequently, when $r\to\infty$,
\begin{equation}\label{approx rate-g}
g(r)\asymp \ln r.
\end{equation}

\noindent\textbf{Proof of Theorem \ref{thm: log}: }
(1)
The main idea is similar to that of the proof of Theorem \ref{thm: p moment convergence rate}. We will use Lemma \ref{le:equi CC'} and different truncating functions to analyze the convergency of $\widetilde{C}_t(\gamma)$. First, for the continuous part, by the Burkholder-Davis-Gundy inequality,
\begin{align}\label{finite for c}
&{\mathbb P}_\mu\left[\left(\sup_{t>1}\int_1^{t}e^{-\lambda s}s^\gamma\int_E \phi(x)S^C(ds, dx)\right)^2\right]\\
&\lesssim
\sup_{t>1}\mathbb P_\mu\int_1^t e^{-2\lambda s} s^{2\gamma}ds\int_E\alpha(x)\phi(x)^2 {X}_s(dx)\nonumber\\
&\lesssim  \int_1^\infty e^{-\lambda s}s^{2\gamma}ds\int_E\alpha(x)\phi(x)^2\nu(dx)\nonumber<\infty.
\end{align}
 For the jump part, we still handle the `small jumps' and the `large jumps'
separately. Define
$$
N^{(1)}:=\sum_{0<\triangle \overline X_s(\phi)<e^{\lambda s}s^{-\gamma}}\delta_{(s, \triangle {\overline X}_s)}\quad \mbox{and}\quad
N^{(2)}:=\sum_{\triangle \overline X_s(\phi)\ge e^{\lambda s}s^{-\gamma}}\delta_{(s, \triangle{\overline X}_s)},
$$
and denote the compensators of $N^{(1)}$  and $N^{(2)}$ by $\widehat N^{(1)}$ and $\widehat N^{(2)}$ respectively.
We write $S^{(J, 1)}$ and $S^{(J,2)}$ for the corresponding martingale measures.  For the `large jumps',
\begin{align*}
&\mathbb P_\mu\Big|\int_1^{\infty}e^{-\lambda s}s^\gamma\int_E \phi(x) S^{(J,2)}(ds, dx)\Big|\\
&\leq  2\mathbb P_\mu\int_1^\infty e^{-\lambda s}s^\gamma  ds\int_E {X}_s(dx)
\int_{e^{\lambda s}s^{-\gamma}}^{\infty}r\pi^\phi(x,dr)\\
&\lesssim   \int_1^\infty s^\gamma ds\int_E\nu(dx)\int_{e^{\lambda s}s^{-\gamma}}^\infty r\pi^\phi(x,dr)\\
&=
\int_1^{R\vee 1} s^\gamma ds\int_E\nu(dx)\int_{e^{\lambda s}s^{-\gamma}}^\infty r\pi^\phi(x,dr)
+\int_{R\vee 1}^\infty s^\gamma ds\int_E\nu(dx)\int_{e^{\lambda s}s^{-\gamma}}^\infty r\pi^\phi(x,dr)\\
&:=I+II,
\end{align*}
where $R>\gamma/\lambda$ is a number such that \eqref{eq: R} holds for $r>R$.  It is obvious that $I<\infty$, we only need to investigate the finiteness of $II$.
Recall that $f$ is defined by \eqref{def-f} and $g$ is the inverse of $f$ on $(\gamma/\lambda, \infty)$.
Applying Fubini's Theorem,
\begin{eqnarray*}
II\le
\int_E\nu(dx)\int_{f(\gamma/\lambda)}^\infty r\pi^\phi(x,dr) \int_{0}^{g(r)} s^\gamma ds
\end{eqnarray*}
It follows from \eqref{approx rate-g}  that
\[
\int_{0}^{g(r)} s^\gamma ds=\dfrac{g(r)^{\gamma+1}}{\gamma+1}\asymp (\ln r)^{\gamma+1} \qquad \mbox{for}\,\,  r>R.
\]
Thus when \eqref{assum: log moment} holds, we have
\[
\sup_{t>1} \mathbb P_\mu\Big|\int_1^{t}e^{-\lambda s}s^\gamma\int_E \phi(x) S^{(J,2)}(ds, dx)\Big|<\infty.
\]
Therefore the process $\int_1^{t}e^{-\lambda s}s^\gamma\int_E \phi(x) S^{(J,2)}(ds, dx)$ converges  $\mathbb P_\mu${-a.s.} and in $L^1(\mathbb P_\mu)$.
Now let us analyze the `small jumps' part.
\begin{align*}
& {\mathbb P}_\mu\left[\Big(\sup_{t>1}\int_1^{t}e^{-\lambda s}s^\gamma\int_E \phi(x) S^{(J,1)}(ds, dx)\Big)^2\right]\\
&=\int_{1}^{\infty}e^{-2\lambda s} s^{2\gamma} ds
\int_E P_s^\beta\Big(\int_0^{f(s)} r^2\pi^\phi(\cdot,dr)\Big)(y)\mu(dy)\\
&\lesssim \int_{1}^{\infty}s^{\gamma}/f(s) ds\int_E \nu(dx)\int_0^{f(s)} r^2\pi^\phi(x,dr)\\
&\leq \int_{1}^{\infty}s^{\gamma}/f(s) ds\int_E \nu(dx)\int_0^{1} r^2\pi^\phi(x,dr)
+\int_{1}^{1\vee R}s^{\gamma}/f(s) ds\int_E \nu(dx)\int_1^{1\vee f(s)} r^2\pi^\phi(x,dr)\\
&\quad +\int_{1\vee R}^\infty s^{\gamma}/f(s) ds\int_E \nu(dx)\int_1^{1\vee f(s)} r^2\pi^\phi(x,dr)\\
&:=III+IV+V.
\end{align*}
It is easy to check that both $III$ and $IV$ are finite. Applying Fubini's theorem in $V$, we get
\[
V\leq \int_E \nu(dx)\int_1^{\infty} r^2\pi^\phi(x,dr)\int_{g(r)}^\infty s^{\gamma}/f(s) ds.
\]
Let $H(r)=\int_{g(r)}^\infty s^{\gamma}/f(s) ds$, then $\lim_{r\rightarrow\infty} H(r)=0$. Note that as $r\to\infty$,
\[
H'(r)=\dfrac{g(r)^\gamma g'(r)}{f(g(r))}\asymp \dfrac{(\ln r)^\gamma}{r^2}.
\]
Thus $H(r)\asymp \dfrac{(\ln r)^\gamma}{r}$ as $r\to\infty$. Therefore, $V<\infty$ when \eqref{assum: log moment} holds.   Hence it follows that
the martingale $ \int_1^{t}e^{-\lambda s}s^{\gamma}\int_E \phi(x) S^{(J,1)}(ds, dx)$ converges $\mathbb P_\mu$-a.s. and in
$L^2(\mathbb P_\mu)$ as $t\to\infty$.
In conclusion, when the moment condition \eqref{assum: log moment} holds, the martingale $\widetilde{C}_t(\gamma)$ converges $\mathbb P_\mu$-almost surely and in $L^1(\mathbb P_\mu)$ as $t\to\infty$.  It follows from Lemma \ref{le:equi CC'} that
\[
 \int_0^t s^{\gamma-1}(M_{\infty}(\phi)-M_s(\phi))ds \quad{\rm converges}\quad \mathbb P_\mu\mbox{-a.s.}
\]
and $ M_\infty(\phi)-M_t(\phi)=o(t^{-\gamma})$, $\mathbb P_\mu$-a.s. as $t\to\infty$.
In particular, when $\gamma\geq 1$,
\[
\int_0^\infty (M_\infty(\phi)-M_t(\phi)) dt<\infty, \qquad \mathbb P_\mu\mbox{-a.s.}
\]

(2) Now let us consider the case that
$\int_E\nu(dx)\int_1^\infty r(\ln r)^{1+\gamma}\pi^{\phi}(x,dr)=\infty$.
Without loss of generality, we may assume that
\begin{equation}\label{assm:log}
\int_E\nu(dx)\int_1^\infty r(\ln r)^{\gamma}\pi^{\phi}(x,dr)<\infty.
\end{equation}
In fact, if
\begin{equation*}
\int_E\nu(dx)\int_1^\infty r(\ln r)^{\gamma}\pi^{\phi}(x,dr)=\infty,
\end{equation*}
then by assumption \eqref{LlogL},
$\gamma>1$. Therefore there is some $\tilde\gamma>0$ and some integer $n>0$ such that $\gamma=n+\tilde \gamma$,
\begin{equation}\label{assm:log'}
\int_E\nu(dx)\int_1^\infty r(\ln r)^{1+\tilde\gamma}\pi^{\phi}(x,dr)=\infty
\end{equation}
and
\begin{equation*}
\int_E\nu(dx)\int_1^\infty r(\ln r)^{\tilde\gamma}\pi^{\phi}(x,dr)<\infty.
\end{equation*}
If we can prove that $C_t(\widetilde\gamma)$ does not converge  as $t\to\infty$,
then $C_t(\gamma)$ does not converge either.

Let $\widehat{N}^{(2)}$ be the compensator of $N^{(2)}$. Then for any non-negative Borel function $F$ on
$\mathbb R_+\times\mathcal M(E_\partial)$,
\[
\int_0^\infty \int_{\mathcal{M}(E_\partial)} F(s,\upsilon)\widehat N^{(2)}(ds, d\upsilon)
=\int_0^\infty ds\int_E {X}_s(dx)\int_{f(s)}^\infty F(s,r\phi(x)^{-1}\delta_x)\pi^\phi(x,dr).
\]
Define a measure $L(ds,dx)$ on $[0,\infty)\times E$ such that for any non-negative Borel function $g$ on $\mathbb R_+\times E$,
\[
 \int_0^\infty\int_Eg(s,x)L(ds,dx)=\int_0^\infty \int_{\mathcal{M}(E_\partial)} F_g(s,\upsilon)\widehat N^{(2)}(ds, d\upsilon),
\]
which is equivalent to
\[
\int_0^\infty\int_Eg(s,x)L(ds,dx)=\int_0^\infty ds\int_E\phi^{-1}(x){X}_s(dx) \int^\infty_{f(s)}rg(s, x)\pi^\phi(x,dr).
\]
Suppose $\mu\in \mathcal{M}^0(E)$.
We claim that as  $t\to\infty$,
\begin{equation}\label{claim-K-limit}
 K_t(\gamma):=\int_1^t e^{-\lambda s}s^{\gamma}\int_E\phi(x)\left[M(ds,dx)+L(ds,dx)\right] \mbox{ converges }
 \mathbb P_{\mu}\mbox{-a.s.}
\end{equation}
 In  fact, for the continuous part of $M$, by \eqref{finite for c},
$\int_1^t e^{-\lambda s}s^{\gamma}\int_E\phi(x)S^C(ds,dx)$ converges
$\mathbb P_{\mu}$-almost surely as $t\to\infty$.
For the `small jump' part, using  the arguments for the `small jumps' in (1), assumption \eqref{assm:log} is enough to guarantee that $\int_1^t e^{-\lambda s}s^{\gamma}\int_E\phi(x)S^{(J,1)}(ds,dx)$ converges
$\mathbb P_{\mu}$-almost surely as $t\to\infty$.
We are left to analyze the `big jumps' part. Thanks to assumption \eqref{assm:log},
\begin{align*}
&\mathbb P_\mu\left(\sum_{\stackrel{\Delta \overline X_s(\phi)\geq f(s)}{s>1}}1\right)
=\mathbb P_\mu\left(\int_1^\infty  ds\int_E{X}_s(dx)\int_{f(s)}^\infty \pi^\phi(x,dr)\right)\\
&=\int_E\mu(dy)\int_1^\infty ds P^\beta_s\Big( \int_{f(s)}^\infty \pi^{\phi}(\cdot,dr)\Big)(y)\\
&\lesssim  \mu(\phi)\int_E\nu(dx)\int_1^\infty e^{\lambda s} ds \int_{f(s)}^\infty \pi^\phi(x,dr)\\
&\leq  \int_E\nu(dx) \int_{f(\frac{\gamma}{\lambda})}^\infty \pi^\phi(x,dr)\int_0^{g(r)}e^{\lambda s}ds\\
&\lesssim  \int_E\nu(dx)\int_{1}^\infty r(\ln^+ r)^\gamma \pi^\phi(x,dr)
<\infty,
\end{align*}
where the second to last inequality comes from \eqref{approx rate-g} and the fact that
$$
\int_0^{g(r)}e^{\lambda s}ds=\dfrac{1}{\lambda}f(s)s^\gamma\Big|_0^{g(r)}=\dfrac{1}{\lambda}(rg(r)^\gamma-1).
$$
Thus the measure $N^{(2)}$ is a finite measure.  Consequently we have as $t\to\infty$,
\begin{equation}\label{large-int-finite}
\int^t_1\int_{\mathcal{M}(E)}F_{e^{-\lambda s}s^\gamma\phi(x)}(s,\upsilon)N^{(2)}(ds, d\upsilon)\to
\sum_{\stackrel{\Delta \overline X_s(\phi)\geq (s)}{s>1}}{}e^{-\lambda s}s^\gamma\Delta \overline X_s(\phi)<\infty,
\end{equation}
since the sum is a finite sum.  Now \eqref{large-int-finite} implies our claim
\eqref {claim-K-limit}.

 Set $L_t=\int_0^te^{-\lambda s}\int_E\phi(x)L(ds,dx)$
and let $L_\infty$ denote its increasing limit.
Then
$$
L_\infty=\int_0^\infty e^{-\lambda s}\int_E\phi(x)L(ds,dx).
$$
We first claim that $L_\infty<\infty, \, \mathbb P_\mu$-a.s. In fact, by the definition  \eqref{def-f} of $f$, $f(s)\geq f(\gamma/\lambda)$ for any $s>0$. Thus it follows from \eqref{assum: moment finite} that for any $s>0$,
\[
\int_{f(s)}^\infty r\pi^{\phi}(x,dr)\leq \int_{f(\gamma/\lambda)}^\infty r\pi^{\phi}(x,dr)\lesssim \phi(x).
\]
Thus
\begin{align*}
 &\int_0^{\gamma/\lambda} e^{-\lambda s}\int_E\phi(x)L(ds,dx)
=\int_0^{\gamma/\lambda} e^{-\lambda s}ds\int_E \overline X_{s-}(dx)\int_{f(s)}^\infty r\pi^{\phi}(x,dr)\\
&\lesssim  \int_0^{\gamma/\lambda} e^{-\lambda s}ds\int_E \phi(x)X_{s}(dx)=\int_0^{\gamma/\lambda}M_s(\phi)ds<\infty,\quad \mathbb P_\mu\mbox{-a.s.}
\end{align*}
By Assumption \ref{asp:H2},
\begin{align*}
 &\mathbb P_\mu\left(\int_{\gamma/\lambda}^\infty e^{-\lambda s}\int_E\phi(x)L(ds,dx)\right)=\mathbb P_\mu\left(\int_{\gamma/\lambda}^\infty e^{-\lambda
s}ds\int_E \overline X_{s-}(dx)\int_{f(s)}^\infty r\pi^{\phi}(x,dr)\right)\\
&=\int_{\gamma/\lambda}^\infty e^{-\lambda s}ds\int_E \mu(dy)P^\beta_s\left(\int_{f(s)}^\infty r\pi^{\phi}(\cdot,dr)\right)(y)\\
&\lesssim  \int_{\gamma/\lambda}^\infty ds\int_E \nu(dx)\int_{f(s)}^\infty r\pi^{\phi}(x,dr)\\
&=\int_E\nu(dx)\int_{f(\gamma/\lambda)}^\infty r\pi^{\phi}(x,dr)\int_{\gamma/\lambda}^{g(r)} ds\\
&\leq \int_E\nu(dx)\int_{f(\lambda/\gamma)}^\infty r g(r)\pi^{\phi}(x,dr)<\infty,
\end{align*}
which implies our claim.

 Now using Lemma \ref{le:equi CC'} and  Remark \ref{remark for general CC'}, \eqref {claim-K-limit} implies that
$$
\int_0^ts^{\gamma-1} \left(M_\infty(\phi)-M_s(\phi)+L_\infty-L_s\right)ds
$$
converges and $(M_\infty(\phi)-M_t(\phi))+(L_\infty-L_t)=o(t^{-\gamma})$
$\mathbb P_{\mu}$-a.s.
Thus the $\mathbb P_\mu$-almost sure convergence of
$\int_0^ts^{\gamma-1} \left(M_\infty(\phi)-M_s(\phi)\right)ds$ as $t\to\infty$
is equivalent to that $\int_0^ts^{\gamma-1} \left(L_\infty-L_s\right)ds$ converges $\mathbb P_\mu$-almost surely to a finite random variable as
$t\to\infty$, and $M_\infty(\phi)-M_t(\phi)=o(t^{-\gamma})$ if and only if $(L_\infty-L_t)$ does.
Since the integrand is non-negative,
we always have limit $\int_1^\infty s^{\gamma-1} \left(L_\infty-L_s\right)ds\leq\infty$.
If we can prove that,  under the assumption \eqref{assum: lower},
\begin{equation}\label{L-limt-infty}
\mathbb{P}_\mu\left(\int_1^\infty s^{\gamma-1} \left(L_\infty-L_s\right)ds=\infty\right)>0.
\end{equation}
Then $\int_0^ts^{\gamma-1} \left(M_\infty(\phi)-M_s(\phi)\right)ds$ does not converge $\mathbb P_\mu$-almost surely as $t\to\infty$.
  Now we are left to prove \eqref{L-limt-infty}.
Note that by \eqref{assum: lower}, there exists $T_2>\max(T_1, \gamma/\lambda)$ such that for $t\ge T_2$,
 \begin{align*}
&L_\infty-L_t=\int_t^\infty e^{-\lambda s}ds\int_E \overline X_{s-}(dx)\int_{f(s)}^\infty r\pi^\phi(x,dr)\\
&\gtrsim\int_t^\infty e^{-\lambda s}ds\int_F\phi(x) \overline X_{s-}(dx)\int_{f(s)}^\infty r\int_E\nu(dy)\pi^\phi(y,dr).
 \end{align*}
Put $\rho(dr)=\int_E\nu(dy)\pi^\phi(y,dr)$. Then for $t\ge T_2$,
\begin{align}\label{est:Linf-L}
&L_\infty-L_t\gtrsim \int_t^\infty e^{-\lambda s}X_{s}(\phi 1_F)ds\int_{f(s)}^\infty r \rho(dr)\\
&\gtrsim \sum_{n=1+[t]}^\infty \int_n^{n+1} ds e^{-\lambda s}X_s(\phi 1_F)\int_{f(s)}^\infty r \rho(dr)
\nonumber\\
&\gtrsim \sum_{n=1+[t]}^\infty \int_{f(n+1)}^\infty r \rho(dr)\int_n^{n+1} e^{-\lambda s}X_s(\phi 1_F)ds.
\nonumber
\end{align}
The third inequality comes from the fact that $f(s)$ is an increasing function for $s>\gamma/\lambda$. It is easy to check that
$\int_0^{1+T_2}t^{\gamma-1}dt
\sum_{n=1+[t]}^\infty \int_{f(n+1)}^\infty r \rho(dr)\int_n^{n+1} e^{-\lambda s}X_s(\phi I_F)ds<\infty$ $\mathbb P_\mu$-almost surely. So
$\int_{1+T_2}^\infty t^{\gamma-1}dt
\sum_{n=1+[t]}^\infty \int_{f(n+1)}^\infty r \rho(dr)\int_n^{n+1} e^{-\lambda s}X_s(\phi 1_F)ds$ have the same convergence property as $\int_0^\infty t^{\gamma-1}dt\sum_{n=1+[t]}^\infty \int_{f(n+1)}^\infty r \rho(dr)\int_n^{n+1} e^{-\lambda s}X_s(\phi 1_F)ds$. By Theorem
\ref{thm intconv} in the Appendix,
the two integrals above converge almost surely on $\{M_\infty(\phi)>0\}$ if and only if the following integral
$$
\int_0^\infty t^{\gamma-1}dt\sum_{n=1+[t]}^\infty \int_{f(n+1)}^\infty r \rho(dr)
$$
is finite. Exchanging the order of integration, we obtain
\begin{align*}
&\int_{0}^\infty t^{\gamma-1}dt\sum_{n=1+[t]}^\infty \int_{f(n+1)}^\infty r \rho(dr)\geq \int_{0}^\infty t^{\gamma-1}dt\int_{t+1}^\infty ds \int_{f(s)}^\infty r \rho(dr)\\
&=\int_1^\infty ds\int_0^{s-1}dt\int_{f(s)}^\infty r\rho(dr)=\dfrac{1}{\gamma}\int_1^\infty(s-1)^{\gamma}ds\int_{f(s)}^\infty r\rho(dr)\\
&\geq\dfrac{1}{\gamma}\int_{f(\gamma/\lambda)}^\infty r\int_{(\gamma/\lambda)\vee 1}^{g(r)\vee 1} (s-1)^\gamma ds=\dfrac{1}{\gamma(\gamma+1)}\int_{f(\gamma/\lambda)}^\infty r[(g(r)-1)^{\gamma+1}-(\gamma/\lambda\vee 1-1)^{\gamma+1}]\rho(dr)\\
&=\infty.
\end{align*}
The last equality is due to that $g(r)\asymp \ln r$ as $r\to\infty$ and \eqref{assum: log moment infty}.
Thus on the event $\{M_\infty(\phi)>0\}$, which has positive probability,
\[
\int_{0}^\infty s^{\gamma-1} \left(L_\infty-L_s\right)ds=\infty,
\]
 almost surely. Thus \eqref{L-limt-infty} is valid and the proof is complete.

 Now we analyze the convergence rate of $L_\infty-L_t$.
It follows from \eqref{est:Linf-L}, Theorem \ref{thm intconv} and the monotonicity of $f$ that,
on $\{M_\infty(\phi)>0\}$, almost surely
 when $t\ge T_2$,
 \[
L_\infty-L_{t}
 \gtrsim \int_{t}^\infty ds \int_{f(s)}^\infty r \rho(dr)=\int_{f(t)}^\infty r[g(r)-t]\rho(dr).
 \]
If $L_\infty-L_t=o(t^{-\gamma})$, then $\int_{f(t)}^\infty r[g(r)-t]\rho(dr)=o(t^{-\gamma})$, or equivalently,
$\int_{t}^\infty r[g(r)-g(t)]\rho(dr)=o(g(t)^{-\gamma})$ as $t\to\infty$.
From \eqref{eq: R},
\[
\int_{t}^\infty r[g(r)-g(t)]\rho(dr)=\int_{t}^\infty r\rho(dr)\int_t^rg'(u)du\asymp \int_{t}^\infty r[\ln r-\ln t]\rho(dr).
\]
By \eqref{approx rate-g},
$\int_{t}^\infty r[g(r)-g(t)]\rho(dr)=o(g(t)^{-\gamma})$ is equivalent to
$\int_{t}^\infty r[\ln r-\ln t]\rho(dr)=o((\ln t)^{-\gamma})$.
 Conversely, when $\int_{t}^\infty r[\ln r-\ln t]\rho(dr)=o((\ln t)^{-\gamma})$
as $t\to\infty$ does not hold,
$L_\infty-L_t=o(t^{-\gamma})$ does not hold almost surely.
Consequently $M_\infty(\phi)-M_t(\phi)=o(t^{-\gamma})$ does not hold almost surely.
\qed

\section{Appendix }
In this appendix, we prove the following result used in the proof of Theorem
\ref{thm: log}.

\begin{thm}\label{thm intconv}
For any Borel subset $F$ of $E$ and $\mu\in\mathcal M(E)$,
\[
\lim_{n\to\infty}\int_n^{n+1} e^{-\lambda s}\langle\phi 1_F,X_s\rangle ds=\langle\phi 1_F,\nu\rangle M_\infty(\phi),\quad  \mathbb P_\mu\mbox{-a.s.}
\]
\end{thm}

The proof of this theorem is based on the following five results.  The idea of the proof is mainly from \cite{LRS1}.
For any $n\in\mathbb N, u>0$, and $h\in\mathcal B_b^+(E)$,
define
$$
H_{n+u}(h):=e^{-\lambda (n+u)}
\int_0^{n+u}\int_EP^{\beta}_{(n+u)-s}(\phi h)(x)
S^{(1,1)}(ds,dx),
$$
$$
L_{n+u}(h):= e^{-\lambda (n+u)}\int_0^{n+u}\int_E P^{\beta}_{(n+u)-s}(\phi h)(x)
S^{(2,1)}(ds, dx),
$$
and
$$
C_{n+u}(h):=e^{-\lambda (n+u)}
\int_0^{n+u}\int_E(P^\beta_{(n+u)-s}\phi h)(x) S^{C}(ds,dx).
$$
\begin{lemma}\label{limit lemma 1}
If $\int_El(x)\nu(dx)<\infty$, then for any $u>0$, $\mu\in {\mathcal M}(E)$ and
$h\in\mathcal B_b^+(E)$,
\begin{equation}\label{finite 2}
\sum_{n=1}^\infty \mathbb P_\mu\left[H_{n+u}(h) -\mathbb
P_\mu(H_{n+u}(h)\big|\mathcal F_{n})\right]^2<\infty
\end{equation}
and
\begin{equation}\label{limit1}
\lim_{n\rightarrow\infty} \left(H_{n+u}(h)- \mathbb P_\mu[H_{n+u}(h)\big|\mathcal F_{n}]\right)=0,\quad \mbox{ in } L^2(\mathbb
P_\mu)  \mbox{ and } \mathbb P_\mu\mbox{-a.s.}
\end{equation}
Moreover,
\begin{equation}\label{limit1ofint}
\lim_{n\rightarrow\infty} \int_0^1\left(H_{n+u}(h)- \mathbb P_\mu[H_{n+u}(h)\big|\mathcal F_{n}]\right) du=0,\quad  \mathbb P_\mu\mbox{-a.s.}
\end{equation}
\end{lemma}
\noindent {\bf Proof:} Since $P^\beta_t(\phi h)$ is  bounded on $[0,T]\times E$ for any $T>0$, the process
$$
H_{t}(h):=e^{-\lambda (n+u)} \int_0^{t} \int_E P^{\beta}_{(n+u)-s}(\phi h)(x) S^{(1,1)}(ds,dx),\qquad
t\in [0, n+u]
$$
is a martingale with respect to $(\mathcal F_{t})_{t\leq n+u}$. Thus
$$
\mathbb P_\mu(H_{n+u}(h)\big|\mathcal F_n)=e^{-\lambda (n+u)} \int_0^{n}\int_E P^{\beta}_{(n+u)-s}(\phi h)(x) S^{(1,1)}(ds,dx),
$$
and hence
\begin{equation}\label{martingale-difference1}
H_{n+u}(h)-\mathbb P_\mu(H_{n+u}(h) \big|\mathcal
F_n) =e^{-\lambda(n+u)}\int_{n}^{n+u}\int_E P^{\beta}_{(n+u)-s}(\phi h)(x) S^{(1,1)}(ds, dx).
\end{equation}
Since
\begin{align*}
&M_t:=e^{-\lambda (n+u)}\int_{n}^t\int_EP^{\beta}_{(n+u)-s}(\phi h)(x) S^{(1,1)}(ds,dx)\\
&=\int_{n}^t \int_{\mathcal M(E_\partial)}F_{e^{-\lambda(n+u)}P^{\beta}_{(n+u)-\cdot}(\phi h)}(s,\upsilon) (N^{(1,1)}-
\widehat{N}^{(1,1)})(ds,d\upsilon), \quad t\in[n, n+u]
\end{align*}
is a martingale with quadratic variation
\begin{eqnarray*}
\int_{n}^t\int_{\mathcal M_f(E_\partial)} F^2_{e^{-\lambda(n+u)}P^{\beta}_{(n+u)-\cdot}(\phi h)}(s,v)
\widehat{N}^{(1,1)}(ds,dv),
\end{eqnarray*}
we have
\begin{align}\label{inequality 1}
&\mathbb P_\mu\left[\left[H_{n+u}(h)-\mathbb P_\mu(H_{n+u}(h)\big|\mathcal F_{n})\right]^2\right]\\
&= \mathbb P_\mu\int_{n}^{n+u}\int_{\mathcal M(E_\partial)}F^2_{e^{-\lambda (n+u)}P^{\beta}_{(n+u)-\cdot}(\phi h)}(s,\upsilon) \widehat{N}^{(1,1)}(ds,d\upsilon)\nonumber\\
 &=\mathbb P_\mu\bigg[\sum_{s\in \widetilde J^{(1,1)}_{n,u}}
F^2_{e^{-\lambda (n+u)}P^{\beta}_{(n+u)-\cdot}(\phi h)}(s, \Delta \overline X_s)\bigg],
\nonumber
\end{align}
where $\widetilde J^{(1,1)}_{n,u}= J^{(1,1)}\bigcap\ [n,n+u]$.  Note that for any $h\in\mathcal B_b^+(E)$,
\begin{equation}\label{uniform boundary2}
P^\beta_t(\phi h)(y)\leq \|h\|_{\infty}e^{\lambda t}\phi(y), \quad \forall\ t\ge 0,\  y\in E.
\end{equation}
Therefore we obtain
\begin{align*}
&\mathbb P_\mu\bigg[\sum_{s\in \widetilde
J^{(1, 1)}_{n, u}
}F^2_{e^{-\lambda(n+u)} P^{\beta}_{(n+u)-\cdot}(\phi h)}(s, \Delta \overline X_s)\bigg]\\
&=\mathbb P_{\mu}\int_{n}^{n+u} ds\int_E {X}_s(dx)
\int_0^{e^{\lambda s}}F^2_{e^{-\lambda (n+u)}P^{\beta}_{(n+u)-\cdot}(\phi h)}(s, r\phi(x)^{-1}\delta_x)\pi^\phi(x,dr)\\
&=e^{-2\lambda (n+u)}\int_{n}^{n+u}ds\int_E \mu(dy)P^{\beta}_s\left(\int_0^{e^{\lambda s}} [P^\beta_{(n+u)-s} (\phi h)(\cdot)\phi(\cdot)^{-1}]^2r^2\pi^\phi(\cdot,dr)\right)(y)\\
&\le \|h\|_\infty^2 \int_{n}^{n+u}e^{-2\lambda s}ds \int_E \mu(dy)P^{\beta}_s\left(\int_0^{e^{\lambda s}} r^2\pi^\phi(\cdot,dr)\right)(y),
\end{align*}
where in the second equality we used the fact that
\begin{eqnarray}\label{formu-F}
F_{e^{-\lambda (n+u)}P^{\beta}_{(n+u)-\cdot}(\phi h)}(s,r\phi(x)^{-1}\delta_x)= re^{-\lambda(n+u)}\phi^{-1}(x)P^{\beta}_{(n+u)-s}(\phi h)(x)
\end{eqnarray}
and in the last  inequality we used \eqref{uniform boundary2}. From the fact that $\int_0^{e^{\lambda s}} r^2\pi^\phi(\cdot,dr)$ is integrable with respect to $\nu$ for any $s>0$, it follows that
\begin{align*}
&\mathbb P_\mu\bigg[\sum_{s\in \widetilde J^{(1, 1)}_{n,u}}
F^2_{e^{-\lambda(n+u)}P^\beta_{(n+u)-\cdot}(\phi h)}(s, \Delta \overline X_s)\bigg]\\
&\leq C \|h\|_\infty^2\langle\phi,\mu\rangle\int_E\nu(dx)\int_{n}^\infty e^{-\lambda s}ds\int_0^{e^{\lambda s}}r^2 \pi^\phi(x, dr).
\end{align*}
Summing over $n$, we get
\begin{align}\label{I+II}
&\sum_{n=1}^\infty \mathbb P_\mu\bigg[\sum_{s\in \widetilde J^{(1,1)}_{n,u}}
F_{e^{-\lambda (n+u)}P^\beta_{(n+u)-\cdot}(\phi h)}(s,\Delta \overline X_s)^2\bigg]\\
&\leq \sum_{n=1}^\infty  C\| h\|_\infty^2\langle \phi,\mu\rangle\int_E\nu(dx)\int_{n}^\infty e^{-\lambda s}ds\int_0^{e^{\lambda s}}r^2\pi^\phi(x,dr)\nonumber\\
&\leq  C\|h\|_\infty^2\langle\phi,\mu\rangle\int_E\nu(dx)\int_0^\infty dt\int_{t}^\infty  e^{-\lambda s}ds\int_0^{e^{\lambda s}}r^2 \pi^\phi(x, dr)\nonumber\\
&= C\|h\|_\infty^2\langle\phi,\mu\rangle\int_E\nu(dx)\int_{0}^\infty s e^{-\lambda s}ds\int_0^{e^{\lambda s}}r^2 \pi^\phi(x,dr)\nonumber\\
&\le  C\|h\|_\infty^2\langle \phi,\mu\rangle\int_E\nu(dx)\int_1^\infty r^2 \pi^\phi(x,dr)\int_{\lambda^{-1}\ln r}^\infty se^{-\lambda s}ds\nonumber\\
&+C\|h\|_\infty^2\langle \phi,\mu\rangle\int_E\nu(dx)\int_0^1 r^2 \pi^\phi(x,dr)\int_{ 0}^\infty se^{-\lambda s}ds\nonumber\\
&=:I+II.\nonumber
\end{align}
Using the assumption $\int_E\nu(dx)\int_0^\infty (r\wedge r^2) \pi^\phi(x,dr)<\infty$, we immediately get that
$II<\infty$. On the other hand,
$$
I=\frac{C}{\lambda^2}\|h\|_\infty^2\langle
\phi,\mu\rangle\int_E\nu(dx)\int_1^\infty r(\ln r+1) \pi^\phi(x,dr).
$$
Now we can use $\int_E l(x)\nu(dx)<\infty$ and $\int_E\nu(dx)\int_0^\infty (r\wedge r^2) \pi^\phi(x,dr)<\infty$ again to get that $I<\infty$. The proof of
\eqref{finite 2} is now complete.   For any $\varepsilon>0$,  using
\eqref{finite 2} and Chebyshev's inequality we have
\begin{align*}
&\sum_{n=1}^\infty\mathbb P_\mu\left(\big|H_{n+u}(h)- \mathbb
P_\mu[H_{n+u}(h)\big|\mathcal F_n]\big|>\varepsilon\right)\\
&\leq \varepsilon^{-2}\sum_{n=1}^\infty \mathbb P_\mu\left[H_{n+u}(h) -\mathbb P_\mu(H_{n+u}(h)\big|\mathcal F_n)\right]^2
< \infty.
\end{align*}
Then \eqref{limit1} follows easily from the  Borel-Cantelli lemma. \eqref{limit1ofint} follows
from \eqref{I+II} and the bounded convergence theorem.
\qed

\begin{lemma}\label{limit lemma 1(2)}
If $\int_El(x)\nu(dx)<\infty$, then for any $u>0$, $\mu\in  {\mathcal M}(E)$ and
$h\in\mathcal B_b^+(E)$ we have
\begin{equation}\label{limit2}
\lim_{n\rightarrow\infty}L_{n+u}(h)-\mathbb P_\mu\left[L_{n+u}(h)\big|\mathcal F_{n}\right]=0,\qquad \mbox{ in } L^1(\mathbb P_\mu)\,  \mbox{ and }\, \mathbb P_\mu\mbox{-a.s.}
\end{equation}
and
\begin{equation}\label{limit2ofint}
\lim_{n\rightarrow\infty}\int_0^1\left(L_{n+u}(h)-\mathbb P_\mu\left[
L_{n+u}(h)\big|\mathcal F_{n}\right]\right)du=0,\quad  \mathbb P_\mu\mbox{-a.s.}
\end{equation}
\end{lemma}

\noindent{\bf Proof.}  It is easy to see that
$$
\mathbb P_\mu\left[ L_{n+u}(h)\big|\mathcal F_{n}\right]=e^{-\lambda(n+u)} \int_0^{n}\int_EP^{\beta}_{(n+u)-s}(\phi h)(x) S^{(2,1)}(ds,dx).
$$
Therefore,
\begin{align}\label{identity 3}
&\left|L_{n+u}(h)-\mathbb P_\mu\left[L_{n+u}(h)\big|\mathcal F_{n}\right]\right|\\
&=\left| e^{-\lambda(n+u)} \int_{n}^{n+u}\int_EP^\beta_{(n+u)-s}(\phi h)(x)S^{(2,1)}(ds,dx)\right|\nonumber\\
&=\left|\int_{n}^{n+u}\int_{\mathcal M(E_\partial)}F_{e^{-\lambda(n+u)}P^{\beta}_{(n+u)-\cdot}(\phi h)(\cdot)}(s, \upsilon)(N^{(2,1)}-\widehat{N}^{(2,1)})( ds, d\upsilon)\right|\nonumber\\
&\leq  \int_{n}^{n+u}\int_{\mathcal M(E_\partial)}F_{e^{-\lambda(n+u)}P^{\beta}_{(n+u)-\cdot}(\phi h)(\cdot)}(s,\upsilon)(N^{(2,1)}+\widehat{N}^{(2,1)})(ds,d\upsilon).\nonumber
\end{align}
Using \eqref{uniform boundary2} we get,
\begin{align*}
&\int_{n}^{n+u}\int_{\mathcal M(E_\partial)}F_{e^{-\lambda (n+u)}P^{\beta}_{(n+u)-\cdot}(\phi h)(\cdot)}(s,\upsilon)(N^{(2,1)}+\widehat{N}^{(2,1)})(ds,d\upsilon)\\
&\leq \int_{n}^\infty\int_{\mathcal M(E_\partial)}F_{\|h\|_\infty e^{-\lambda \cdot}\phi}(s,\upsilon)(N^{(2,1)}+\widehat{N}^{(2,1)})(ds,d\upsilon)
\end{align*}
Taking expectation, we obtain
\begin{align*}
&\mathbb P_\mu\int_{n}^\infty\int_{\mathcal M(E_\partial)}F_{\|h\|_\infty e^{-\lambda \cdot}\phi}(s, \upsilon)(N^{(2,1)}+\widehat{N}^{(2,1)})(ds,d\upsilon)\\
&=2\|h\|_\infty\mathbb P_\mu \left[\int_{n}^\infty e^{-\lambda s}ds\int_E {X}_s(dx)
\int_{e^{\lambda s}}^\infty  r \pi^{\phi}(x,dr)\right]\\
&=2\|h\|_\infty\int_{n}^\infty e^{-\lambda s}ds \int_E\mu(dy)P^\beta_s\left(\int_{e^{\lambda s}}^\infty r\pi^\phi(\cdot,dr)\right)(y)\\
&\leq 2C\|h\|_\infty\langle\phi,\mu\rangle\int_{n}^\infty ds\int_E\nu(dx)\int_{e^{\lambda s}}^\infty r \pi^{\phi}(x, dr)\\
&\leq  2C\|h\|_\infty\langle\phi,\mu\rangle\int_E\nu(dx)\int_{e^{\lambda n}}^\infty r \pi^{\phi}(x,dr)\int_0^{\lambda^{-1}\ln r}ds\\
&= \frac{2C\|h\|_\infty}{\lambda} \langle\phi,\mu\rangle\int_E\nu(dx)\int_{e^{\lambda n}}^\infty  r\ln r \pi^{\phi}(x,dr).
\end{align*}
Note that $\int_E\nu(dx)\int_{e^{\lambda n}}^\infty  r\ln r \pi^{\phi}(x,dr)\le\int_El(x)\nu(dx)<\infty$.
Applying the dominated convergence theorem and using the fact that
$$
\int_{n}^\infty\int_{\mathcal M(E_\partial)}F_{\|h\|_\infty e^{-\lambda \cdot}\phi}(s, \upsilon)(N^{(2,1)}+\widehat{N}^{(2,1)})(ds,d\upsilon)
$$
is decreasing in $n$, we obtain that, when $\int_El(x)\nu(dx)<\infty$,
\begin{equation}\label{limit-L}
\lim_{n\rightarrow\infty}\int_{n}^\infty\int_{\mathcal M(E_\partial)}
F_{\|h\|_\infty e^{-\lambda \cdot}\phi}(s, \upsilon)(N^{(2,1)}+\widehat{N}^{(2,1)})(ds,d\upsilon) =0,\quad \mbox{ in } L^1(\mathbb
P_\mu)  \mbox{ and } \mathbb P_\mu\mbox{-a.s.}
\end{equation}
Therefore by \eqref{identity 3}, we have \eqref{limit2} and \eqref{limit2ofint}.
The proof is complete.\qed

\begin{lemma}\label{cont-conv}
 For any $u>0$,
$\mu\in  {\mathcal M}(E)$ and
$h\in\mathcal B_b^+(E)$ we have
\begin{equation}\label{limit5}
\lim_{n\rightarrow\infty}C_{n+u}(h)-\mathbb P_\mu[C_{n+u}(h)
\big|\mathcal F_{n}]=0,\quad \mbox{ in } L^2(\mathbb
P_\mu)  \mbox{ and } \mathbb P_\mu\mbox{-a.s.}
\end{equation}
and
\begin{equation}\label{limit5ofint}
\lim_{n\rightarrow\infty}\int_0^1\left(C_{n+u}(h)-\mathbb P_\mu[C_{n+u}(h)
\big|\mathcal F_{n}]\right)du=0,\quad  \mathbb P_\mu\mbox{-a.s.}
\end{equation}
\end{lemma}
\noindent {\bf Proof:}
Note  that
$$
\mathbb P_\mu(C_{n+u}(h)\big|\mathcal
F_{n})=e^{-\lambda(n+u)}
\int_0^{n}\int_E(P^\beta_{(n+u)-s}\phi h)(x) S^{C}(ds,dx),
$$
and then
\begin{equation}\label{martingale-difference3}
 C_{n+u}(h)-\mathbb P_\mu[C_{n+u}(h)\big|\mathcal F_{n}]=  e^{-\lambda(n+u)}\int_{n}^{n+u}\int_E(P^{\beta}_{(n+u)-s}\phi h)(x) S^{C}(ds,dx).
\end{equation}
From the quadratic variation formula, it follows that
\begin{align}
&\mathbb P_\mu\left[\left[C_{n+u}(h)-\mathbb P_\mu(C_{n+u}(h)
\big|\mathcal F_{n})\right]^2\right]\\
&=\int_{n}^{n+u}e^{-2\lambda(n+u)}ds\int_E
\mu(dx)P_s^\beta\left(\left(P^\beta_{(n+u)-s}\phi h\right)^2(\cdot)\alpha(\cdot)\right)(x) \nonumber\\
&\le\|h\|^2_\infty\int_{n}^{n+u}e^{-\lambda s}ds
\int_EP^{\beta}_s(\alpha\phi^2)(x)\mu(dx)\nonumber\\
&\leq  \frac{1}{\lambda}\|\phi\alpha\|_{\infty}\|h\|^2_\infty\langle
\phi,\mu\rangle e^{-\lambda n}.\nonumber
\end{align}
Therefore, we have
\begin{equation}
\sum_{n=1}^\infty \mathbb P_\mu\left[C_{n+u}(h)-\mathbb
P_\mu(C_{n+u}(h)\big|\mathcal F_{n})\right]^2<\infty.
\end{equation}
By the Borel-Cantelli lemma, we get \eqref{limit5}.
\eqref{limit5ofint} follows from
the bounded convergence theorem.
\qed

\begin{thm}\label{main skelton}
If $\int_El(x)\nu(dx)<\infty$, then for any $\mu\in {\mathcal M}(E)$ and $h\in\mathcal B_b^+(E)$
we have
$$
\lim_{n\rightarrow\infty}e^{-\lambda n}\langle\phi h,X_{n}\rangle=M_{\infty}(\phi)\int_E\phi(z)h(z)\nu(dz) ,\quad\mbox{in } L^1(\mathbb P_\mu)\mbox{ and } \mathbb P_\mu\mbox{-a.s.}
$$
\end{thm}

\noindent{\bf Proof:}
By combining the three lemmas above, we can easily get that, if $\int_El(x)\nu(dx)<\infty$, then for any
$u>0$,
$\mu\in {\mathcal M}(E)$ and
$h\in\mathcal B_b^+(E)$ we have
\begin{equation}\label{limitmain1}
\lim_{n\rightarrow\infty} e^{-\lambda(n+u)}\langle\phi h,
X_{n+u}\rangle -\mathbb P_\mu\left[
e^{-\lambda(n+u)}\langle\phi h,X_{n+u}\rangle\big|\mathcal F_{n}\right]=0,\quad
\mbox{ in } L^1(\mathbb P_\mu)  \mbox{ and } \mathbb P_\mu\mbox{-a.s.}
\end{equation}
Indeed, $e^{-\lambda(n+u)}\langle\phi h, X_{n+u}\rangle$ can
be decomposed into four parts:
\begin{align*}
&e^{-\lambda(n+u)}\langle\phi h, X_{n+u}\rangle\\
&=e^{-\lambda (n+u)}\langle P^{\beta}_{n+u}(\phi h),\mu\rangle+e^{-\lambda(n+u)}\int_0^{n+u}\int_EP^{\beta}_{(n+u)-s}(\phi h)(x)M(ds,dx)\\
&=e^{-\lambda (n+u)}\langle P^\beta_{n+u}(\phi h),\mu\rangle+H_{n+u}(h)+L_{n+u}(h)+ C_{n+u}(h).
\end{align*}
Therefore,
\begin{align*}
&e^{-\lambda(n+u)}\langle\phi h, X_{n+u}\rangle-\mathbb P_\mu\left[ e^{-\lambda(n+u)}\langle\phi h, X_{n+u}\rangle\big|\mathcal F_{n}\right]\\
&=  H_{n+u}(h)-\mathbb P_\mu\left[H_{n+u}(h)|\mathcal F_{n}\right]+L_{n+u}(h)-\mathbb P_\mu\left[L_{n+u}(h)|\mathcal F_{n}\right]\\
  &\quad+C_{n+u}(h)-\mathbb P_\mu\left[C_{n+u}(h)|\mathcal F_{n}\right].
\end{align*}
Now \eqref{limitmain1} follows immediately from
Lemmas \ref{limit lemma 1}--\ref{cont-conv}.
By the mean formula and the Markov property of
superprocesses, we have
\begin{equation}\label{cond-id}
\mathbb P_\mu\left[ e^{-\lambda(n+u)}\langle\phi h,
X_{n+u}\rangle\big|\mathcal F_{n}\right]
=e^{-\lambda n}\langle e^{-\lambda u}
P^{\beta}_{u}(\phi h), X_{n}\rangle.
\end{equation}
By Assumption 2,
\[
e^{-\lambda u}P^{\beta}_{u}(\phi h)(x)=\phi(x)\langle\phi h, \nu\rangle(1+C_{u,x,\phi h}),
\]
where $\left|C_{u,x,\phi h}\right|\leq c_u$ and $\lim_{u\to\infty}c_u=0$.  Hence,
\begin{eqnarray}
e^{-\lambda n}\langle e^{-\lambda u}P^{\beta}_{u}(\phi
h), X_{n}\rangle &\geq& (1-c_u) e^{-\lambda n}\langle\phi, X_{n}\rangle \langle\phi h, \nu\rangle\nonumber\\
&=& (1-c_u)M_{n}(\phi)\langle\phi h, \nu\rangle,\quad \mathbb P_\mu\mbox{-a.s.}\label{rseqn1}
\end{eqnarray}
and
\begin{eqnarray}
e^{-\lambda n}\langle e^{-\lambda m}P^{\beta}_{m}(\phi h), X_{n}\rangle &\leq& (1+c_m) e^{-\lambda n}\langle\phi, X_{n}\rangle \langle\phi h, \nu\rangle\nonumber\\
&=& (1+c_m)M_{n}(\phi)\langle\phi h, \nu\rangle,\quad \mathbb P_\mu\mbox{-a.s.}\label{rseqn2}
\end{eqnarray}
Using \eqref{cond-id},
%Lemma \ref{limit lemma1}
\eqref{limitmain1}
and \eqref{rseqn2}, we get that
for any $u>0$,
\begin{align*}
&\limsup_{n\rightarrow\infty}e^{-\lambda n}\langle\phi h,
X_{n}\rangle=\limsup_{n\rightarrow\infty}e^{-\lambda (n+u)}\langle\phi h, X_{n+u}\rangle\\
&=\limsup_{n\rightarrow\infty}e^{-\lambda n}\langle e^{-\lambda u}P^\beta_{u}(\phi h), X_{n}\rangle
\le \limsup_{n\rightarrow\infty}(1+c_u) M_{n}(\phi)\langle\phi h, \nu\rangle\\
&=(1+c_u)M_\infty(\phi)\langle\phi h, \nu\rangle.
\end{align*}
Letting $u\to\infty$, we get
\begin{equation}\label{rseqn3}
\limsup_{n\rightarrow\infty}e^{-\lambda n}\langle\phi h,
X_{n}\rangle\le M_\infty(\phi)\langle\phi h, \nu\rangle.
\end{equation}
Similarly, using \eqref{cond-id},
%Lemma \ref{limit lemma1}
\eqref{limitmain1}
and \eqref{rseqn1} we
get that for any $u>0$,
\begin{align*}
&\liminf_{n\rightarrow\infty}e^{-\lambda n}\langle\phi h,
X_{n}\rangle=\liminf_{n\rightarrow\infty}e^{-\lambda (n+u)}\langle\phi h,
X_{n+u}\rangle\\
&=\liminf_{n\rightarrow\infty}e^{-\lambda n}\langle e^{-\lambda u}
P^\beta_{u}(\phi h), X_{n}\rangle
\ge \liminf_{n\rightarrow\infty}(1-c_u)M_{n}(\phi)\langle\phi h, \nu\rangle\\
&=(1-c_u)M_\infty(\phi)\langle\phi h, \nu\rangle.
\end{align*}
Letting $u\to\infty$, we get
\begin{equation}\label{rseqn4}
\liminf_{n\rightarrow\infty}e^{-\lambda n}\langle\phi h, X_{n}\rangle\ge M_\infty(\phi)\langle\phi h, \nu\rangle.
\end{equation}
Combining \eqref{rseqn3} and \eqref{rseqn4} we arrive at the almost sure assertion of the theorem.
Since $e^{-\lambda n}\langle\phi h, X_{n}\rangle$ is controlled by
a constant multiple of $M_{n}(\phi)$, which is uniformly integrable by Proposition \ref{degeneracy theorem},
the $L^1$ assertion now follows immediately from the almost sure assertion.
\qed

\noindent{\bf Proof of Theorem \ref{thm intconv}:} For any $s>n$
\begin{align*}
&e^{-\lambda s}\langle\phi 1_F,X_s\rangle \\
&=e^{-\lambda s}\langle P^{\beta}_{s-n}(\phi 1_F),X_n\rangle+e^{-\lambda s}\int_n^{s}P^{\beta}_{s-u}(\phi 1_F)(x)M(du,dx)\\
&=e^{-\lambda s}\langle P^{\beta}_{s-n}(\phi 1_F),X_n\rangle+\left(H_{s}(\phi 1_F)- \mathbb P_\mu[H_{s}(\phi 1_F)\big|\mathcal F_{n}]\right)+
\left(L_{s}(\phi 1_F)-\mathbb P_\mu\left[L_{s}(\phi 1_F)|\mathcal F_{n}\right]\right)\\
&\quad+C_{s}(\phi 1_F)-\mathbb P_\mu\left[C_{s}(\phi 1_F)|\mathcal F_{n}\right].
\end{align*}
Hence,
\begin{align*}
&\int_n^{n+1}e^{-\lambda s}\langle\phi 1_F,X_s\rangle ds\\
&=\int_n^{n+1}e^{-\lambda s}\langle P^{\beta}_{s-n}(\phi 1_F),X_n\rangle ds+\int_n^{n+1}\left(H_{s}(\phi 1_F)- \mathbb P_\mu[H_{s}(\phi 1_F)\big|\mathcal F_{n}]\right)ds \\
&\quad+\int_n^{n+1}\left(L_{s}(\phi 1_F)-\mathbb P_\mu\left[L_{s}(\phi 1_F)|\mathcal F_{n}\right]\right)ds
+\int_n^{n+1}C_{s}(\phi 1_F)-\mathbb P_\mu\left[C_{s}(\phi 1_F)|\mathcal F_{n}\right] ds\\
&= e^{-\lambda n}\langle \left(\int_0^{1}e^{-\lambda s}P^{\beta}_{s}(\phi 1_F) ds\right),X_n\rangle+\int_0^{1}\left(H_{n+s}(\phi 1_F)- \mathbb P_\mu[H_{n+s}(\phi 1_F)\big|\mathcal F_{n}]\right)ds\\
&\quad + \int_0^{1}\left(L_{s+n}(\phi 1_F)-\mathbb P_\mu\left[L_{s+n}(\phi 1_F)|\mathcal F_{n}\right]\right)ds
+\int_0^{1}C_{n+s}(\phi 1_F)-\mathbb P_\mu\left[C_{n+s}(\phi 1_F)|\mathcal F_{n}\right] ds.\\
&=I_n+II_n+III_n+IV_n.
\end{align*}
It has been shown in Lemma \ref{limit lemma 1}, Lemma \ref{limit lemma 1(2)} and Lemma \ref{cont-conv} that
\[
\lim_{n\to\infty}II_n+III_n+IV_n=0.
\]
Since $\int_0^{1}e^{-\lambda s}P^{\beta}_{s}(\phi I_F)(x) ds\leq \phi(x)$, by Theorem \ref{main skelton},
\[
\lim_{n\to\infty}e^{-\lambda n}\langle \left(\int_0^{1}e^{-\lambda s}P^{\beta}_{s}(\phi 1_F) ds\right),X_n\rangle=M_\infty(\phi)
\langle\int_0^{1}e^{-\lambda s}P^{\beta}_{s}(\phi 1_F)(x) ds, \nu\rangle=M_\infty(\phi)\langle\phi I_F,\nu\rangle.
\]
\qed

\medskip

%{\bf Acknowledgements.}

\vspace{.1in}

\begin{singlespace} 
\end{singlespace}
\end{doublespace}

\vskip 0.3truein \vskip 0.3truein

 \noindent {\bf Rong-Li Liu:}  Mathematics and Applied Mathematics, Beijing Jiaotong University, Beijing 100044, P. R. China.
 E-mail: {\tt rlliu@bjtu.edu.cn}
  \\

\bigskip
 \noindent {\bf Yan-Xia Ren:} LMAM School of Mathematical Sciences \& Center for
Statistical Science, Peking University, Beijing 100871, P. R. China.
   E-mail: {\tt yxren@math.pku.edu.cn} \\

\bigskip
\noindent {\bf Renming Song:} Department of Mathematics,
 University of Illinois,  Urbana, IL 61801 U.S.A..
 E-mail: {\tt rsong@illinois.edu} \\
 \bigskip
\end{document}